\documentclass[a4paper, 12pt,reqno]{amsart}
\usepackage[left=3cm,right=3cm]{geometry}

\usepackage[T1]{fontenc} 
\usepackage{ae,aecompl} 
\usepackage{lmodern}

\usepackage[english]{babel}

\usepackage{hyperref}

\usepackage{dynkin-diagrams}

\usepackage{amsthm, amssymb, amsfonts, amsmath, mathrsfs}
\usepackage[all]{xy}

\usepackage{caption}
\usepackage{url}
\usepackage{tikz}

\usepackage{xcolor}

\usepackage{array}
\newcommand{\PreserveBackslash}[1]{\let\temp=\\#1\let\\=\temp}
\newcolumntype{C}[1]{>{\PreserveBackslash\centering}p{#1}}
\newcolumntype{R}[1]{>{\PreserveBackslash\raggedleft}p{#1}}
\newcolumntype{L}[1]{>{\PreserveBackslash\raggedright}p{#1}}

\newcounter{stepcounter}
\theoremstyle{plain}

\newtheorem{thm}{Theorem}[section]

\newtheorem{lemma}[thm]{Lemma}
\newtheorem{prop}[thm]{Proposition}
\newtheorem{cor}[thm]{Corollary}
\newtheorem{conj}[thm]{Conjecture}

\theoremstyle{definition}

\newtheorem{example}[thm]{Example}
\newtheorem{defn}[thm]{Definition}

\newtheorem{remark}[thm]{Remark}

\date{}

\newcommand\bit{\begin{itemize}}
\newcommand\eit{\end{itemize}}
\newcommand\bet{\begin{enumerate}}
\newcommand\eet{\end{enumerate}}
\newcommand\ed{\end{document}}

\DeclareFontFamily{U}{mathx}{\hyphenchar\font45}
\DeclareFontShape{U}{mathx}{m}{n}{
      <5> <6> <7> <8> <9> <10>
      <10.95> <12> <14.4> <17.28> <20.74> <24.88>
      mathx10
      }{}
\DeclareSymbolFont{mathx}{U}{mathx}{m}{n}
\DeclareFontSubstitution{U}{mathx}{m}{n}
\DeclareMathAccent{\widecheck}{0}{mathx}{"71}
\DeclareMathAccent{\wideparen}{0}{mathx}{"75}

\def\cs#1{\texttt{\char`\\#1}}

\newcommand{\e}{\varepsilon}
\newcommand\Om{\Omega}

\newcommand\DEL{\Delta}
\newcommand\LL{\Lambda}

\newcommand\bC{{\mathbb C}}
\newcommand\bN{{\mathbb N}}

\newcommand\F{{\mathcal F}}

\newcommand\fg{{\mathfrak g}}

\newcommand\fh{{\mathfrak h}}

\DeclareMathOperator{\Ad}{Ad}

\DeclareMathOperator{\can}{can}

\newcommand\co{\mathrm{co}}

\newcommand\exd{\mathrm{d}}

\newcommand\End{\mathrm{End}}

\newcommand\unit{\mathrm{U}}
\newcommand\counit{\mathrm{C}}

\newcommand\id{\mathrm{id}}
\newcommand\proj{\mathrm{proj}}

\newcommand\vol{\mathrm{vol}}

\newcommand\ahol{^{(0,1)}}

\newcommand\inv{^{-1}}

\newcommand\by{\times}

\newcommand\oby{\otimes}

\newcommand\sseq{\subseteq}
\newcommand\tl{\triangleleft}
\newcommand\tr{\triangleright}

\def\qbinom#1#2{\ensuremath{\left[\kern-.3em\left[\genfrac{}{}{0pt}{}{#1}{#2}\right]\kern-.3em\right]_q}}

\usepackage{tikz}
\usetikzlibrary{decorations.pathreplacing}

\newcommand{\rGr}{\mathrm{Gr}}

\newcommand{\cO}{\mathcal{O}}
\newcommand{\cA}{\mathcal{A}}

\newcommand{\Zee}{\mathbb{Z}}
\newcommand{\Cee}{\mathbb{C}}
\newcommand{\Ree}{\mathbb{R}}

\newcommand{\eps}{\varepsilon}
\newcommand{\fB}{\mathfrak{B}}
\newcommand{\fsl}{\mathfrak{sl}}
\newcommand{\fsp}{\mathfrak{sp}}
\newcommand{\fo}{\mathfrak{o}}
\newcommand{\fgl}{\mathfrak{gl}}
\newcommand{\fe}{\mathfrak{e}}
\newcommand{\fl}{\mathfrak{l}}
\newcommand{\fk}{\mathfrak{k}}

\newcommand{\cT}{\mathcal{T}}

\newcommand{\ltr}{\triangleleft}

\DeclareMathOperator{\Rep}{\mathsf{Rep}}
\newcommand{\Mod}{\mathsf{Mod}}
\DeclareMathOperator{\YD}{\mathsf{YD}}

\def\clap#1{\hbox to 0pt{\hss#1\hss}}
\def\mathllap{\mathpalette\mathllapinternal}
\def\mathrlap{\mathpalette\mathrlapinternal}

\def\mathllapinternal#1#2{%
  \llap{$\mathsurround=0pt#1{#2}$}}
\def\mathrlapinternal#1#2{%
  \rlap{$\mathsurround=0pt#1{#2}$}}

\newsavebox\qModFirst
\newsavebox\qModSecond
\newcommand{\modz}[2]{
  \sbox{\qModFirst}{$#1$}
  \sbox{\qModSecond}{$#2$}
  \ifdim\wd\qModFirst>\wd\qModSecond
    {}^{\phantom{#1}\mathllap{#1}}_{\phantom{#1}\mathllap{#2}}%
  \else
    {}^{\phantom{#2}\mathllap{#1}}_{\phantom{#2}\mathllap{#2}}%
  \fi\mathsf{mod}_0}
\newcommand{\Modz}[2]{
  \sbox{\qModFirst}{$#1$}
  \sbox{\qModSecond}{$#2$}
  \ifdim\wd\qModFirst>\wd\qModSecond
    {}^{\phantom{#1}\mathllap{#1}}_{\phantom{#1}\mathllap{#2}}%
  \else
    {}^{\phantom{#2}\mathllap{#1}}_{\phantom{#2}\mathllap{#2}}%
  \fi\mathsf{Mod}_0}
\newcommand{\qMod}[4]{{%
      \sbox{\qModFirst}{$#1$}
      \sbox{\qModSecond}{$#2$}
      \ifdim\wd\qModFirst>\wd\qModSecond
        {}^{\phantom{#1}\mathllap{#1}}_{\phantom{#1}\mathllap{#2}}%
      \else
        {}^{\phantom{#2}\mathllap{#1}}_{\phantom{#2}\mathllap{#2}}%
      \fi
  \Mod^{#3}_{#4}}}
\newcommand{\lMod}[2]{\qMod{#1}{#2}{}{}}
\newcommand{\qmod}[4]{{%
      \sbox{\qModFirst}{$#1$}
      \sbox{\qModSecond}{$#2$}
      \ifdim\wd\qModFirst>\wd\qModSecond
        {}^{\phantom{#1}\mathllap{#1}}_{\phantom{#1}\mathllap{#2}}%
      \else
        {}^{\phantom{#2}\mathllap{#1}}_{\phantom{#2}\mathllap{#2}}%
      \fi
  \mathsf{mod}^{#3}_{#4}}}

\newcommand{\rYD}[1]{\YD^{#1}_{#1}}

\newcommand{\what}{\widehat}

\newcommand{\barM}{\overline{M}}

\title{Nichols Algebras and Quantum Principal Bundles}

\author[A. O. Krutov]{Andrey O. Krutov}
\address{Institute of Mathematics, Czech Academy of Sciences, \v{Z}itn\'a 25, 115 67 Prague, Czech Republic}

 \email{krutov@math.cas.cz}

\author[R. \'O Buachalla]{R\'eamonn \'O Buachalla}
\address{Mathematical Institute of Charles University, Sokolovsk\'a 83, 186 75 Prague, Czech Republic} 

\email{obuachalla@karlin.mff.cuni.cz}

\author[K. R. Strung]{Karen R. Strung}
\address{Institute of Mathematics, Czech Academy of Sciences, \v{Z}itn\'a 25, 115 67 Prague, Czech Republic}
\email{strung@math.cas.cz}

\keywords{Nichols algebras, quantum groups, quantum flag manifolds, quantum principal bundles, noncommutative
  geometry}

\subjclass[2020]{
  16T20, 
  46L87, 
  81R60, 
  81R50, 
  17B37, 
  16T05}  

\thanks{
  {\tiny AK was supported through the program ``Oberwolfach Leibniz Fellows'' by the Mathematisches
    Forschungsinstitut Oberwolfach in 2018, by WCMCS post-doctoral fellowship, and by the QuantiXLie Centre of
    Excellence, a project cofinanced by the Croatian Government and European Union through the European
    Regional Development Fund---the Competitiveness and Cohesion Operational Programme (KK.01.1.1.01.0004).
    R\'OB acknowledges FNRS support through  a postdoctoral fellowship within the framework of the MIS Grant
    ``Antipode'' grant number F.4502.18. R\'OB is currently supported by supported by the Charles University PRIMUS grant \emph{Spectral Noncommutative Geometry of Quantum Flag Manifolds} PRIMUS/21/SCI/026.  KRS was supported by a Radboud Excellence Initiative postdoctoral fellowship and Sonata 9 NCN grant 2015/17/D/ST1/02529.
    AK and KRS are supported by the GA\v{C}R project 20-17488Y and \mbox{RVO: 67985840}.}
  }

\begin{document}

\maketitle

\begin{abstract}
We introduce a general framework for associating to a homogeneous quantum principal bundle a Yetter--Drinfeld module structure on the cotangent space of the base calculus. The holomorphic and anti-holomorphic Heckenberger--Kolb calculi of the quantum Grassmannians are then presented in this framework. This allows us to express the calculi in terms of the corresponding Nichols algebras. The extension of this result to all irreducible quantum flag manifolds is then conjectured.
\end{abstract}

\section{Introduction}


Exterior and symmetric algebras play a fundamental and ubiquitous role in classical differential
geometry. Their quantum counterparts, however, are much more poorly understood, as is their role in
noncommutative differential geometry. The most developed framework we have for understanding exterior and symmetric algebras in the noncommutative setting is the theory of Nichols algebras, an important class of braided Hopf algebras. Nichols algebras first appeared in~\cite{Nichols78} as a tool for the constructing new examples of Hopf algebras. Later they would arise independently in a number of works, for example, the work of Woronowicz on differential calculi~\cite{WoroDC,WoronowiczDC89} and the work of Majid on braided groups~\cite{Maj91,Maj93}. Nichols algebras were subsequently used  to give an abstract construction of quantised enveloping algebras~\cite{Lusztig1993book,MajLeabh}.  Most famously, Nichols algebras are basic invariants of pointed Hopf algebras, and are crucial in Andruskiewitsch and Schneider's remarkable classification program for Hopf algebras~\cite{AndruskiewitschSchneider1998,AndruskiewitschSchneiderAnnals,AndruskiewitschICM}. Other notable applications include the work of Bazlov~\cite{Bazlov2006}, where Nichols algebras were used to describe the cohomology rings of the flag manifold $G/B$ of a semisimple Lie group $G$. This was later generalised to the case of the (small) quantum cohomology ring for $G/B$ in~\cite{KirillovMaeno2005}. Nichols algebras have also seen applications in conformal field theory~\cite{Semikhatov2012,SemikhatovTipunin2012,SemikhatovTipunin2013} and more recently~\cite{Lentner2021}.


Woronowicz rediscovered Nichols algebras as part of his investigation of bicovariant differential calculi over
Hopf algebras. He showed that the left invariant forms of a bicovariant calculus possess a Yetter--Drinfeld
structure, and hence the structure of a braided vector space. He then took the associated Nichols algebra,
endowed it with a commutator differential, producing an extension of the first-order calculus to a
differential graded algebra.
 Such an approach leads to noncommutative analogues of curvature, Bianchi identities, and Lie derivatives~\cite{AschieriCastellani1993}.
However, an obvious problem is that any bicovariant calculus over a~quantum coordinate algebra $\cO_q(G)$ will not have classical dimension, even for the simplest example of $\cO_q(SU_2)$~\cite{SchSchulerClass}. As a result applications, of Nichols algebras to the noncommutative geometry of quantum groups remain limited.


In contrast to the situation for Drinfeld--Jimbo quantum groups, the framework of covariant differential
calculi has proved itself to be an ideal  setting for investigating the noncommutative geometry of
quantum flag manifolds.
In particular, the seminal work of Heckenberger and Kolb has shown that the irreducible quantum flag manifolds admit an essentially unique covariant $q$-deformation of their classical de Rham complex \cite{HK, HKdR}. These differential calculi arguably constitute the most important family of noncommutative differential structures in the theory of quantum groups.  In the special case of the $A$-series, the irreducible quantum flag manifolds are precisely the quantum Grassmannians, which properly contain the quantum projective spaces $\cO_q(\mathbb{CP}^n)$. In particular, when $n = 1$, this includes the celebrated Podle\'s sphere $\mathcal{O}_q(\mathbb{CP}^1) \cong \mathcal{O}_q(S^2)$.

Quantum flag manifolds are constructed as subspaces of right coinvariant elements with respect to a
quantum Levi subgroup $\cO_q(L_S)$. As such they only admit a left $\cO_q(G)$-coaction, meaning that for these
quantum spaces one can only speak of left $\cO_q(G)$-covariant calculi.  This places the Heckenberger--Kolb
calculi outside Woronowicz's bicovariant framework. However, left covariance of the calculi gives their
cotangent spaces the structure of an $\cO_q(L_S)$-comodule. Since $\cO_q(L_S)$ is a coquasitriangular Hopf algebra, its category of comodules is a braided
monoidal category with braiding~$\sigma_R$.
Therefore any $\cO_q(L_S)$-comodule~$V$ will have an associated Nichols algebra~$\fB(V,\sigma_R)$,
 $q$-deforming the classical exterior algebra. Unfortunately, in
most cases the dimension of~$\fB(V,\sigma_R)$ will not be the same as the dimension of the classical exterior algebra,
for example, $\fB(V,\sigma_R)$ will often be infinite-dimensional. To address this non-classical behaviour, Berenstein and Zwicknagl
introduced the novel notion of \emph{quantum} exterior algebras~\cite{BerensteinZwicknagl2008}. Their approach is
set in the framework of coboundary categories (as opposed to braided monoidal categories) and cactus groups
(as opposed to braid groups). This solves the problem of non-classical dimension for a distinguished family of
$U_q(\fg)$-modules classified in~\cite{Zwicknagl2009}. In particular, for the irreducible quantum flag manifolds, Berenstein and Zwicknagl's quantum
exterior algebras have classical dimension. These quantum exterior
algebras would later be used by Kr\"ahmer and Tucker-Simmons in their construction of noncommutative
Dolbeault--Dirac operators over the irreducible quantum flag manifolds~\cite{MTSUK}.


In this paper, we return to the more standard braid group approach, producing for the first time a Nichols
algebra presentation of the Heckenberger--Kolb calculi $\Omega^1_q(\mathrm{Gr}_{n,m})$ of the quantum
Grassmannians. We introduce a general approach, based around the theory  of quantum principal bundles, for
producing a Yetter--Drinfeld braiding on the cotangent space of a covariant first-order calculus over a
quantum homogeneous space. Indeed, our approach can be considered as a generalisation of Woronowicz's work to
the setting of quantum homogeneous spaces. This general approach is then applied to the quantum principal
bundle description of $\Omega^1_q(\mathrm{Gr}_{n,m})$ introduced in~\cite{BWGrass}, which uses a differential
calculus on $\cO_q(SU_n)$ constructed  from the coquasitriangular structure of $\cO_q(SU_n)$, following the
approach of~\cite[\textsection5]{BeggsMajidBook}. For recent advances, see~\cite{AschieriFioresiLatini2021,Branimir2021CMP,BrzezinskiSzymanski2021} and
references therein. The resulting Yetter--Drinfeld structure is shown to be both diagonal and of Hecke type for the special case of quantum projective space, but non-diagonal and of non-Hecke type for the quantum Grassmannians $\cO_q(\text{\rm Gr}_{n,m})$, with $m \notin \{1,n-1\}$. 


Next, we investigate the associated Nichols algebras, starting with the degree two terms where the braiding
gives an explicit description of the  commutation  relations of the calculi. The remaining terms of the Hilbert--Poincar\'e
series  for the associated Nichols algebra are calculated using  Poincar\'e duality for Nichols algebras,
the fact that the maximal prolongation $\Phi(\Omega^{(0,\bullet)})$ is a Frobenius algebra (as observed
in~\cite{MTSUK}), and quantum Howe duality for the quantum Levi subalgebra of the quantum Grassmannians. The
Hilbert--Poincar\'e series of the Nichols algebra and that of  $\Phi(\Omega^{(0,\bullet)})$ coincide, meaning
that the two graded algebras coincide. This constitutes the principal result of the paper, and confirms that
the braiding produced is novel and distinct from the (co)quasitriangular braiding. Moreover, it establishes an important point of contact between Nichols algebras and noncommutative geometry.


This project leads to a number of new questions. The first and most obvious is how to extend our results to
all the irreducible quantum flag manifolds. Since the quantum enveloping algebras of the other series are also
quasitriangular, it seems reasonable to expect that we can adapt the construction of~\cite{BWGrass} to this
more general setting. Moreover, the calculation of the Hilbert--Poincar\'e series is far easier for cases outside the $A$-series.  As we show in~\S\ref{sec:other}, the fact that the relations of a Nichols
algebra of Hecke type are generated in degree two, means that the Hilbert--Poincar\'e series can be simply
concluded from the action of the braiding on degree two forms.
A further problem is to define Lie derivatives, contraction operators, Maurer--\/Cartan forms and the
corresponding Cartan calculus in the manner of~\cite{AschieriCastellani1993}.
Another challenging, but  important, goal is to produce a Nichols algebra description of  the whole Dolbeault double complex of the irreducible quantum flag manifolds.

\subsection{Summary of Results}


The paper is organised as follows: In~\S\ref{sec:Pre} we recall necessary preliminaries about
braided monoidal categories, Yetter--Drinfeld modules, Nichols algebras, and differential calculi over
Hopf algebras and quantum homogeneous spaces, focusing on complex structures, connections, holomorphic structures, principal comodule algebras, and strong principal connections.

In~\S\ref{sec:main} we establish a general quantum principal bundle framework for constructing a Yetter--Drinfeld module structure on the cotangent space of a covariant first-order differential calculus. Along the way we establish some novel categorical equivalences to clarify the underlying  processes at work.

In~\S\ref{sec:Grass} the basic definitions and results of Drinfeld--Jimbo quantum groups are recalled. We then present the definition of a quantum flag manifold, focusing on the special case of the quantum Grassmannians $\cO_q(\mathrm{Gr}_{n,m})$ and their Heckenberger--Kolb calculi $\Omega^1_q(\mathrm{Gr}_{n,m})$. We apply the general results of~\S\ref{sec:main} to the anti-holomorphic part of $\Omega^1_q(\mathrm{Gr}_{n,m})$ and prove the main result of the paper:

\begin{thm}
For any quantum Grassmannian, the anti-holomorphic tangent space $V^{(0,1)}$ of its Heckenberger--Kolb calculus admits the structure of an $\cO_q(L_S)$-Yetter--Drinfeld module, such that the associated Nichols algebra is isomorphic, as an $\cO_q(L_S)$-comodule algebra, to the maximal prolongation $V^{(0, \bullet)}$. An analogous result holds for the holomorphic tangent space $V^{(1,0)}$.
\end{thm}

In \textsection 5 we conjecture the existence of a Yetter--Drinfeld structure for the cotangent spaces of all the irreducible quantum flag manifolds. We also discuss the conjectured isomorphism between the  Nichols algebra and calculus for each simple Lie algebra series, using explicit representation-theoretic calculations.

\subsubsection*{Acknowledgments}
Part of this work was carried when R\'OB and AK visited KRS at the Institute of Mathematics, Astrophysics and
Particle Physics at Radboud University, Nijmegen and we thank the institute for their support.
KRS and R\'OB are also grateful for a visit to Mathematisches Forschungsintitut Oberwolfach in December 2018
where AK was a Leibniz fellow, and to the Mathematics Department at the University of Zagreb during July 2019.
All three authors benefitted from meeting at the conference ``Quantum Flag Manifolds in Prague'' at the Charles
University in September 2019.
We also thank  Vincent Grandjean and Simon Lentner for helpful discussions, as well as the referees for helpful
remarks and comments.
Special thanks also to Pi\k{e}kna
Mery, who made this collaboration possible.

\section{Preliminaries}\label{sec:Pre}

\subsection{Hopf algebras and braided Hopf algebras}

Throughout this paper~$A$ and~$H$ will denote Hopf algebras, and all Hopf algebras are assumed to unital,
with bijective antipode, and defined over the complex numbers. We denote the coproduct, counit, and antipode
of a Hopf algebra by~$\Delta$, $\e$, and $S$, respectively. Throughout we use Sweedler notation, and write
$a^+ := a - \eps(a)1$, for any $a \in A$, and $V^+ := V \cap \ker(\e)$, for~$V$ a subspace of~$A$.
\subsubsection{Braided Hopf algebras}\label{sec:braidedHopf}

A \emph{braiding} on a monoidal category~$\mathsf{C}$ is a natural
isomophism $\sigma$ between functors $-\otimes -$ and $-\otimes^{\text{op}}-$
such that the relevant hexagonal diagrams commute,
see~\cite[\S8.1]{EtingofTensorCat} for details.

A \emph{braided} monoidal category is a pair consisting of a monoidal category and a~braiding.

Let $\mathsf{C}$ be a braided monoidal category with braiding~$\sigma$. To give the tensor product~$A\otimes B$ the structure of an associative algebra in~$\mathsf{C}$, we define a  multiplication by the formula
\begin{equation}\label{eq:BraidMul}
  m_{A\otimes B} := (m_A\otimes m_B)\circ(\id_A\otimes \sigma_{B,A}\otimes\id_B).
\end{equation}

A \emph{braided bialgebra} in~$\mathsf{C}$ is an object~$A$ in~$\mathsf{C}$ endowed with an associative algebra structure in~$\mathsf{C}$
and a coalgebra structure in~$\mathsf{C}$ such that its coproduct~$\Delta$ and counit~$\eps$
are algebra morphisms with respect to the multiplication in~$A\otimes A$ given by~\eqref{eq:BraidMul}. A \emph{braided Hopf algebra} in~$\mathsf{C}$ is a bialgebra in~$\mathsf{C}$ admitting an   antipode which is a morphism in~$\mathsf{C}$. For further details on braided Hopf algebras, we direct the reader to~\cite{EtingofTensorCat} and~\cite{BakalovKirillovJr}.

\subsubsection{Yetter--Drinfeld modules}

An important example of a braided monoidal category is the category of (right) {\em Yetter--Drinfeld modules} $V$
over a Hopf algebra~$H$, which are those right $H$-modules~$V$ with action~$\tl$, and a right $H$-comodule structure
such that
\begin{equation}\label{eq:YDcondition}
  v_{(0)} \tl h_{(1)} \oby v_{(1)}  h_{(2)} = (v \tl h_{(2)})_{(0)} \oby h_{(1)}(v\tl h_{(2)})_{(1)},
  \quad\text{for $h\in H$, $v\in V$.}
\end{equation}
We denote the category of Yetter--Drinfeld modules, endowed with its obvious monoidal structure, by
$\YD^H_H$. A braiding for the category is defined by
\begin{equation}\label{eq:YDbraidingGeneral}
  \sigma: V \oby W \to W \oby V, \qquad v \oby w \mapsto w_{(0)} \oby v \tl w_{(1)},\qquad\text{for $v\in V$, $w\in W$.}
\end{equation}

Note that for any $V\in\YD^H_H$, the tensor algebra~$\cT(V)$ is a braided Hopf algebra in~$\YD^H_H$ with
\[
  \Delta(v) := v \otimes 1 + 1 \otimes v,\qquad
  S(v) := -v, \qquad \eps(v) :=0,\qquad
  \text{for $v\in V$}.
\]

\subsection{Nichols algebras}
For a detailed introduction to Nichols algebras we refer the reader to the
surveys~\cite{Andruskiewitsch2017intro,AndruskiewitschAngiono2017} and~\cite{HeckenbergerSchneiderBook}. Let $\mathbf{B}_n$ denote the \emph{braid group on $n$ strands}, that is, the group generated by $n-1$ elements $ \beta_1, \dots, \beta_{n-1}$ subject to the relations
\begin{align*}
  \beta_i \beta_{i+1} \beta_i &= \beta_{i+1} \beta_i \beta_{i+1}, &{}& 1 \leq i \leq n-2,\\
  \beta_i \beta_ j &= \beta_j \beta_i, &{}& 1 \leq i,j \leq n-2,\  |i-j| \geq 2 .
\end{align*}

When  $V\in\YD^H_H$ is finite-dimensional as a vector space, we obtain
a representation of the braid group on $n$ strands
\[ \rho_n : \mathbf{B}_n \to GL(V^{\otimes n}),\]
given by 
\[ \rho_n(\beta_i) = \id \otimes \cdots \otimes \id \otimes \sigma \otimes \id \otimes \cdots \otimes \id,\]
where $\sigma$ is acting on $V \otimes V$ in position $i$ and $i+1$.

There is a canonical surjective group homomorphism onto the symmetric group~$\mathbf{S}_n$, 
\[
  \varphi_n:\mathbf{B}_n \to \mathbf{S}_n,
\]
which maps $\beta_i$ to the simple transposition $\tau_i = (i, i+1)$.
Let $\ell(g)$ denote the length of an~element~$g\in\mathbf{S}_n$.
The projection~$\varphi_n$ admits a set-theoretic section, called the {\em Matsumoto section}
\[
  s_n:  \mathbf{S}_n \to \mathbf{B}_n,
\]
which is determined by $s_n(\tau_i) = \beta_i$ and 
\(
s_n(\tau_i \tau_{i+1}) = s_n(\tau_i)s_n(\tau_{i+1}), 
\)
for $1 \leq i \leq n$, and $s_n(gf)  = s_n(g)s_n(f)$ if
\(
\ell(gf) = \ell(g)+\ell(f),
\)
for $g,f\in\mathbf{S}_n$.
Note that $s_n$ is \emph{not} a group homomorphism.

The \emph{braided symmetriser} is given by the map
\[
  \mathfrak{S}^\sigma_n(V) := \sum_{g \in \mathbf{S}_n}\rho_n(s_n(g)) : V^{\otimes n} \to V^{\otimes{n}}.
\]
We denote
\[
  \ker\mathfrak{S}^\sigma(V) := \bigoplus_{n\in\mathbb{Z}_{\geq0}} \ker\mathfrak{S}^\sigma_n(V).
\]

\begin{defn}\label{NichDef} The \emph{Nichols algebra} of~$V$ is the braided Hopf algebra in~$\YD^H_H$ defined by
  \[
    \fB(V,\sigma) := \cT(V) \big/  \ker\mathfrak{S}^\sigma(V).
  \]
\end{defn}
 In what follows we will write~$\fB(V)$ when the braiding on $V$ is clear.

Since $\ker\mathfrak{S}^\sigma(V)$ is a~homogeneous ideal of~$\cT(V)$, the Nichols algebra~$\fB(V)$ has a unique $\mathbb{Z}_{\geq0}$-grading 
\[
  \fB(V) \simeq \bigoplus_{n\in\mathbb{Z}_{\geq0}} \fB_n(V),\]
  where $\fB_n(V) := \cT^n(V) \big/   \ker\mathfrak{S}^\sigma_n(V).$

\subsection{Principal comodule algebras}

For a right $H$-comodule $V$ with structure map~$\Delta_R$, we say that an element $v \in V$
is \emph{(right) coinvariant} if $\Delta_R(v) = v \otimes 1$. We denote the subspace of all $H$-coinvariant elements by
$V^{\co(H)}$, and call it the {\em (right) coinvariant subspace} of the coaction.

A right  {\em $H$-comodule algebra} $(P,\Delta_R)$ is a right $H$-como\-dule  which is also an algebra such that the comodule
structure map ${\Delta_R:P \to P \otimes H}$ is an algebra map. We say that $P$ is an {\em $H$-Hopf--Galois extension of} $B := P^{\co(H)}$ if, for $m_P$ the multiplication of
$P$, an isomorphism $P \oby_B P \simeq P \oby H$ is given by
\begin{align*}
\can := (m_P \oby \id) \circ (\id \oby \Delta_R): P \oby_B P \to P \oby H.
\end{align*}

If  the functor  $P \oby_B -: \lMod{}{B} \to \mathsf{Vect}$, from the category of left $B$-modules to the category of vector spaces, preserves and reflects exact  sequences, then we say that $P$ is  {\em faithfully flat} as a right $B$-module. The definition of faithful flatness  for $P$ as a  left $B$-module is analogous.

\begin{defn}
A {\em principal right $H$-comodule algebra}  is a right $H$-comodule algebra  $(P,\Delta_R)$  such that  $P$ is
an $H$-Hopf--Galois extension of $B := P^{\co(H)}$ and $P$ is faithfully flat as a right and left $B$-module.
\end{defn}

\subsection{Quantum homogeneous spaces}

Let $A$ and $H$ be Hopf algebras, and let $\pi : A \to H$ be a surjective Hopf algebra map. A right $H$-coaction,
giving $A$ the structure of a right $H$-comodule algebra, is given by
\[
  \Delta_R : = (\id \otimes \pi) \circ \Delta :A \to A \otimes H.
\]
We call the coinvariant subspace $B := A^{\co(H)}$  of such a coaction a  {\em quantum homogeneous space}.
In this paper we will exclusively consider quantum homogeneous spaces $B = A^{\co(H)}$ for which $A$ is
faithfully flat as a right $B$-module, as it allows us to use Takeuchi's equivalence, see \textsection
\ref{sec:Tak} below. (We note that similar results hold under much weaker assumptions, see~\cite{Skryabin2007}.)
An important fact is that the coproduct of $A$ restricts to a left $A$-coaction
\[
  \Delta_L:B \to A \otimes B,\qquad   b \mapsto b_{(1)}\otimes b_{(2)},
\]  
giving $B$ the structure of a left $A$-comodule algebra.

A {\em strong bicovariant splitting map} is a unital linear map $i:H\to A$ splitting the projection $\pi:A \to H$ such that
\begin{equation} \label{sbsm}
(i \oby \id) \circ \Delta = \Delta_R \circ i, \qquad (\id \oby i) \circ \Delta = \Delta_L \circ i.
\end{equation}
The existence of a bicovariant splitting map implies that the associated quantum homogeneous space gives a principal comodule algebra. For a more detailed discussion of bicovariant splitting maps  see~\cite[\S24]{MajLeabh} and~\cite[\S5]{BeggsMajidBook}.

We say that a Hopf algebra $H$ is \emph{cosemisimple} if its category of comodules $\lMod{H}{}$ is a
semisimple category. Equivalently, $H$ is cosemisimple if it is the direct sum of its simple
subcoalgebras. The following technical lemma follows, for example, from the proof of Lemma~3.6 in~\cite{HolVBs}.
\begin{lemma}\label{cosemiH}  
For $H$ a cosemisimple Hopf algebra, every surjective Hopf algebra map~$\pi:A \to H$ admits a strong
bicovariant splitting map, and hence the corresponding right $H$-comodule algebra  is a principal comodule algebra. 
\end{lemma}

\subsection{Takeuchi's categorical equivalence}\label{sec:Tak}

In this subsection we recall the form of Takeuchi's equivalence~\cite{Tak}, for a quantum homogeneous space $\pi: A \to H$, best suited to the paper.

For any quantum homogeneous space $B = A^{\co(H)}$, we define  $\qMod{A}{B}{}{B}$ to be the category whose  objects are left \mbox{$A$-comodules} \mbox{$\DEL_L:\mathcal{F} \to A \otimes \mathcal{F}$}, endowed with a $B$-bimodule structure such that $\DEL_L(bfc) = \Delta_L(b)\DEL_L(f)\Delta_L(c)$,  for all  $f \in \mathcal{F}, \, b,c \in B$, and whose morphisms  are left $A$-comodule, $B$-bimodule, maps.

Let $\qMod{H}{}{}{B}$ denote the category with objects left $H$-comodules $\Delta_L: V \to H \otimes V$, endowed with a right $B$-module structure such that $\Delta_L(vb) = v_{(-1)}\pi(b_{(1)}) \otimes v_{(0)}b_{(2)}$, for all $v \in V, \, b \in B$, and whose morphisms are left $H$-comodule maps.

Consider the functor 
\begin{align*}
\Phi:\qMod{A}{B}{}{B} \to \qMod{H}{}{}{B}, & & \mathcal{F} \mapsto  \mathcal{F}/B^+\mathcal{F},
\end{align*}
where the left $H$-comodule structure of $\Phi(\mathcal{\F})$ is given by 
$
\Delta_L[f] := \pi(f_{(-1)})\otimes [f_{(0)}],
$
with square brackets denoting the coset of an element in $\Phi(\mathcal{\F})$. In the other direction, we use the cotensor product $\square_{H}$ to define a functor 
\begin{align*}
\Psi: \qMod{H}{}{}{B} \to  \qMod{A}{B}{}{B}, & & V \mapsto A \,\square_{H} V,
\end{align*}
where the left $A$-comodule structure of $\Psi(V)$ is defined on the first tensor factor, the right $B$-module structure is the diagonal one, and if $\gamma$ is a morphism in $\qMod{H}{}{}{}$, then $\Psi(\gamma) := \id \otimes \gamma$. 

An adjoint equivalence of categories between~$\qMod{A}{B}{}{B}$ and~$\qMod{H}{}{}{B}$, which we call
\emph{Takeuchi's equivalence}, is given by the functors $\Phi$ and $\Psi$, the \emph{unit} natural isomorphism
\begin{align*}
\unit^\pi: \F \to \Psi \circ \Phi(\F), & & f \mapsto f_{(-1)} \otimes [f_{(0)}],
\end{align*}
and the \emph{counit} natural transformation
\[
\counit^\pi := (\e \otimes \id): \Phi \circ \Psi(V) \to V.
\]
The \emph{dimension} $\mathrm{dim}(\F)$ of an object $\F \in \qMod{A}{B}{}{B}$ is the vector space dimension of $\Phi(\F)$. 

An important point to note is that for any $B = A^{\co(H)}$ a quantum homogeneous space, Takeuchi's  equivalence implies an isomorphism of Hopf algebras $\phi: H \to A/B^+A$, such that for $\proj:A \to A/B^+A$ the canonical surjection, $\pi = \phi \circ \proj$. This means that  we necessarily have that $\ker(\pi) = B^+A$.

Consider $\Modz{A}{B}$ the full subcategory of $\qMod{A}{B}{}{B}$ consisting of those  objects $\F$ satisfying $B^+\F = \F B^+$. The corresponding full subcategory $\Modz{H}{}$ of $\qMod{H}{}{}{B}$ is given by objects with the trivial right $B$-action. The category $\Modz{A}{B}$ comes equipped with a monoidal structure given by the tensor product $\otimes_B$. Moreover, with respect to the obvious monoidal structure on $\Modz{H}{}$, Takeuchi's equivalence is readily endowed with the structure of a monoidal equivalence (see~\cite[\textsection 4]{MMF2}). 

Finally, we consider $\modz{A}{B}$ the full subcategory of $\Modz{A}{B}$ whose objects are finitely generated as left $B$-modules, and note that it is a~monoidal subcategory of $\Modz{A}{B}$. The corresponding full monoidal subcategory $\modz{H}{}$ of $\lMod{H}{}$ is given by the finite-dimensional left $H$-comodules.

\subsection{The fundamental theorem of two-sided Hopf modules}\label{section:TwoSidedHopfMods}

In this subsection we consider a special case of Takeuchi's equivalence, namely the fundamental theorem of
two-sided Hopf modules. (This equivalence was originally considered in~\cite[Theorem 5.7]{Schauenburg} using a
parallel but equivalent formulation. See also~\cite{Saraccoa19}.)  For a Hopf algebra $A$, the counit
$\e:A \to \mathbb{C}$ is clearly a~surjective Hopf algebra map. The associated quantum homogeneous space is given by
$A = A^{\co(\mathbb{C})}$. In this case, the category $\qMod{A}{B}{}{B}$ specialises to $\qMod{A}{A}{}{A}$, and the category
$\qMod{H}{}{}{B}$ reduces to the category of right $A$-modules $\Mod_A$.  
For this special case, we find it useful to denote the functor $\Phi$ as
\begin{align*}
F: \qMod{A}{A}{}{A}\to \Mod_A, \qquad \F \mapsto \F/A^+\F,
\end{align*}
Moreover, since the cotensor product over $\mathbb{C}$ is just the usual tensor product $\otimes$, we see that the functor $\Psi$ reduces to 
\[
A \otimes -: \Mod_A  \to \qMod{A}{A}{}{A}, \qquad V \mapsto A \otimes V.
\]
Since faithful flatness is trivially satisfied in this case, we have the following consequence of  Takeuchi's equivalence: The \emph{fundamental theorem of two-sided Hopf modules} states that an adjoint equivalence between the
  categories $\qMod{A}{A}{}{A}$ and $\Mod_A$ is given by the functors $F$ and $A \otimes -$,
  and the unit natural isomorphism
\begin{align*}
  \unit:  {}& \F \to A \otimes F(\F),  & & f \mapsto f_{(-1)} \oby [f_{(0)}],
 \end{align*}
 and the counit natural transformation
 \begin{align*}
  \counit:{}& F(A \otimes V) \to V,    & &\left[ a \oby v \right] \mapsto  \e(a)v.
\end{align*} 

\subsection{Differential calculi}
A {\em differential calculus}  $\big(\Omega^\bullet \simeq \bigoplus_{k \in \bN_0} \Omega^k, \exd\big)$ is a differential graded algebra (dg-algebra) which is generated in degree $0$ as a dg-algebra, that is to say, it is generated as an algebra by the elements $a, \exd b$, for $a,b \in \Omega^0$. We call an element $\omega \in \Omega^{\bullet}$ a \emph{form}, and if $\omega \in \Omega^k$, for some $k \in \mathbb{N}$, then $\omega$ is said to be \emph{homogeneous} of degree $|\omega| := k$.  The product of two forms $\omega,\nu \in \Omega^{\bullet}$ is denoted by $\omega \wedge \nu$, unless one of the forms is of degree $0$, whereupon the product is denoted by juxtaposition.  
For a given algebra $B$, a \emph{differential calculus over} $B$ is a differential calculus such that $\Omega^0 =B$.

\subsubsection{Universal differential calculi}

A {\em first-order differential calculus} over an algebra~$B$ is a pair $(\Omega^1(B),\exd)$, where~$\Omega^1(B)$ is a~$B$-bimodule and $\exd: B \to \Omega^1$ is a linear map for which the {\em Leibniz rule} holds
\[
  \exd(ab)=a(\exd b)+(\exd a)b,\qquad \text{for $a,b \in B$},
\]
and for which~$\Omega^1(B)$ is generated as a~left~$B$-module by those elements of the form~$\exd b$, for~$b \in B$. The {\em universal first-order differential calculus} over $B$ is the pair
$(\Omega^1_u(B), \exd_u)$, where~$\Omega^1_u(B)$ is the kernel of the multiplication map $m_B: B \otimes B \to B$ endowed
with the obvious $B$-bimodule structure, and $\exd_u$ is the map defined by
\[
  \exd_u: B \to \Omega^1_u(B), \qquad b \mapsto 1 \oby b - b \oby 1.
\]
By~\cite[Proposition 1.1]{WoroCQPGs}, every first-order differential calculus over $B$ is of the form
$\left(\Omega^1_u(B)/N, \,\proj \circ \exd_u\right)$, where $N$ is a $B$-subbimodule of $\Omega^1_u(B)$, and
we have denoted by $\proj:\Omega^1_u(B) \to \Omega^1_u(B)/N$ the quotient map. This gives a bijective correspondence
between calculi and subbimodules of $\Omega^1_u(B)$.
Moreover, every first-order differential calculus admits an
extension to its {\em maximal prolongation} differential calculus, which is to say, one from which any other extension can be
obtained by quotienting, see for example~\cite[\textsection2.5]{MMF2}.

For $A$ a Hopf algebra and $B$ a left $A$-comodule algebra, we say that a first-order differential calculus $\Omega^1(B)$ over $B$ is   {\em left $A$-covariant} if there exists a (necessarily unique) map $\Delta_L\colon \Omega^1(B) \to A \oby \Omega^1(B)$ satisfying
\[
  \Delta_L(b\exd b') = \Delta_L(b) (\id \oby \exd)\Delta_L(b'), \qquad\text{for $b,b' \in B$}.
\]
Similarly one can define a \emph{right $A$-covariant} first-order differential calculus over
a~right $A$-comodule algebra.

For the special case of~$A$ considered as a~left $A$-comodule algebra over itself, we note that every
left $A$-covariant first-order differential calculus over~$A$ is an object in the category~$\qMod{A}{A}{}{A}$.

\subsubsection{Covariant differential calculi over homogeneous spaces}
In the case that~$B$ is a quantum homogeneous space of the form $B = A^{\co(H)}$, a left $A$-covariant first-order
differential calculus $(\Omega^1(B),\exd)$ is a  natural object in $\qMod{A}{B}{}{B}$.

Let us note that~$B^+$ can be viewed as an object in $\qMod{H}{}{}{B}$. It is easy to see that the map
\begin{equation}\label{eq:xiB}
  \xi_B\colon \Phi(\Omega^1_u(B)) \to B^+,\qquad
  [b_1\exd b_2] \mapsto \eps(b_1)(b_2)^+,\qquad\text{for $b_1,b_2\in B$}
\end{equation}
is an isomorphism in~$\qMod{H}{}{}{B}$. Throughout the paper we will tacitly identify these two objects.

As was shown by Hermisson in~\cite{Hermisson2002} (see also~\cite{Maj}), we can classify left $A$-covariant first-order
differential calculi in terms of subobjects $I_B\subseteq B^+$ in $\qMod{H}{}{}{B}$. In particular, if
$\Omega^1(B) = \Omega^1_u(B)/N$ then the corresponding subobject is given by ${I_B = \xi_B(\Phi(N))}$.
Moreover we have the following commutative diagram
\[
  \xymatrix{
    \Omega^1(B) \ar[rr]^{(\id\otimes\xi_B)\circ \unit^\pi} &\qquad& A\square_H B^+/I_B\\
    B^+ \ar[u]^{\exd}\ar[urr]_{\quad(\id\otimes [-])\circ\Delta_L}
  }
\]
where $[-]\colon B^+ \to B^+/I_B$ is the canonical projection.

\subsubsection{Bicovariant differential calculi over Hopf algebras}
For the special case of a~trivial quantum homogeneous space, Hermisson's classification reduces to Woronowicz's celebrated
theorem classifying left-covariant calculi over a~Hopf algebra $A$~\cite[Theorem 1.5]{WoroDC}. Namely, left-covariant first-order differential calculi over~$A$ are classified by right ideals of $A^+$.

A first-order differential calculus~$\Omega^1(A) = \Omega^1_u(A)/N$ over~$A$ is called \emph{bicovariant} if it is both a right
and a left $A$-covariant differential calculus.
According to~\cite[Theorem~1.8]{WoroDC}, a left-covariant first-order differential calculus over~$A$ is bicovariant if the
corresponding right ideal~$\xi_A(F(N_A))$ of~$A^+$ is invariant with respect to the right adjoint $A$-coaction~$\Ad_R$ given by
\begin{equation}\label{eq:AdA}
  \Ad_R a = a_{(2)}\otimes S(a_{(1)})a_{(3)},\qquad\text{for $a\in A$}.
\end{equation}
Moreover, when $\Omega(A)$ is bicovariant, $F(\Omega(A))$ has the structure of a Yetter--Drinfeld module
over~$A$, see~\cite[\S13 and~\S14]{KSLeabh}.

\subsection{Quantum principal bundles}
For a right $H$-comodule algebra $(P,\Delta_R)$ with $B:=P^{\co(H)}$, it can be shown  that the extension $B \hookrightarrow P$ is Hopf--Galois if and only if the sequence
\begin{align} \label{qpbexactseq}
0 \longrightarrow P\Omega^1_u(B)P {\buildrel \iota \over \longrightarrow} \Omega^1_u(P) {\buildrel {\mathrm{ver}}\over \longrightarrow} P \oby H^+ \longrightarrow 0,
\end{align}
is exact, where $\Omega^1_u(B)$ is the restriction of $\Omega^1_u(P)$ to $B$, $\iota$ is the inclusion map,
$\textrm{ver}$ is the restriction of $\can$ to $\Om_u^1(P)$. It is useful to note that an explicit presentation of the action of $\textrm{ver}$ is given by $\textrm{ver}(a'\exd a) = a'a_{(1)} \otimes \pi(a^+_{(2)})$.
The following definition generalises this sequence to general calculi which are not necessarily universal~\cite[\S5]{BeggsMajidBook}.

\begin{defn} \label{qpb}
  A {\em quantum $H$-principal bundle} is a triple $(P,\Delta_R,\Omega^1(P))$, where
  \begin{enumerate}
  \item $(P,\Delta_R)$ is a right $H$-comodule algebra such that $P$ is a Hopf--Galois extension of $B:=P^{\co(H)}$,  
  \item $\Omega^1(P) \simeq \Omega^1_u(P)/N$ is a~left $H$-covariant first-order differential calculus over~$P$, 
    \[
    \mathrm{ver}(N) = P \oby I,
  \]
   for $I$ some right ideal of~$H^+$ satisfying $\Ad_H I \subseteq I\otimes H$, where $\Ad_H$ is the adjoint coaction of~$H$.
  \end{enumerate}
\end{defn}

Denote  by $\Omega^1(B)$ the restriction of $\Omega^1(P)$ to $B$. Definition~\ref{qpb} implies that 
\begin{align} \label{qpbexactseq}
0 \longrightarrow P\Omega^1(B)P {\buildrel \iota \over \longrightarrow} \Omega^1(P) {\buildrel {\mathrm{ver}}\over \longrightarrow} P \oby (H^+/I )\longrightarrow 0
\end{align}
is a well-defined exact sequence.

\section{Nichols algebras from left $A$-covariant right $H$-covariant calculus}\label{sec:main}
Let $\Omega^1(B)$ be a left $A$-covariant first-order differential calculus over a quantum homogeneous
space~$B$.  In this section we establish sufficient criteria for $\Omega^1(B)$ which allow one to define a (right)
Yetter--Drinfeld module structure on the space of left-coinvariant forms~$\Phi(\Omega^1(B))$.

The objects of the category~$\qMod{A}{A}{H}{A}$ are $A$-bimodules~$\mathcal{F}$ equipped with an $(A,H)$-bicomodule structure $(\Delta_L,\Delta_R)$ such that 
\begin{align*}
\Delta_L(afa') = \Delta(a)\Delta_L(f)\Delta(a'), & & \Delta_R(afa') = \Delta(a)\Delta_R(f)\Delta(a').
\end{align*}
When $H = A$, the category $\qMod{A}{A}{A}{A}$  is  known as the category of \emph{tetramodules}

The objects of the category~$\YD^H_A$ of \emph{relative Yetter--Drinfeld modules}
(see~\cite[\S4]{HeckenbergerSchneider2013}) are right $A$-modules $V$ equipped with a right $H$-comodule
structure $\Delta_R$ satisfying the compatibility condition
\begin{align*}\label{eq:YDAH}
   v _{(0)}\tl a_{(1)}\otimes  v _{(1)}\pi(a_{(2)})
  = ( v \tl a_{(2)})_{(0)}\otimes \pi(a_{(1)}) ( v  \tl a_{(2)})_{(1)}, & & \textrm{ for all }   v \in V, \, a\in A.
\end{align*}

We now consider a generalisation of the well-known equivalence between tetramodules and Yetter--Drinfeld
modules~\cite{Schauenburg}. For any $V  \in \YD^H_A$ we can endow  $A\otimes V$ with the structure of an
object in $\qMod{A}{A}{}{A}$ by taking the left $A$-Hopf module structure of the first tensor factor and the right $A$-module structure given by the action on the tensor product.
Moreover, we endow $A \otimes V$ with a right $H$-comodule structure $\Delta_R\colon A \otimes V \to A \otimes V\otimes H$ as follows
\[
  \Delta_R(a\otimes v ) := a_{(1)}\otimes v _{(0)}\otimes\pi(a_{(2)}) v _{(1)},
  \qquad\text{for $a\in A$, $ v \in V$.}
\]
It is clear that
\[
  \Delta_R(af) = \Delta_R(a)\Delta_R(f),
  \qquad\text{for } a\in A, f \in A \otimes V,
\]
and that $A \otimes V$ is an~$(A,H)$-bicomodule. Moreover, we see that
\begin{align}
  \Delta_R((b\otimes v ) a) ={}
  & \Delta_R(ba_{(1)} \otimes  v \tl a_{(2)})\notag\\
  {}={} &  b_{(1)}a_{(1)}\otimes ( v \tl a_{(3)})_{(0)}\otimes\pi(b_{(2)})\pi(a_{(2)})( v \tl a_{(3)})_{(1)},
  \label{eq:DeltaLHS}\\
  \Delta_R(b\otimes v ) \Delta_R(a) ={}
  & \left(b_{(1)}\otimes v _{(0)}\otimes\pi(b_{(2)}) v _{(1)}\right)
    \left( a_{(1)}\otimes\pi(a_{(2)}) \right)\notag\\
  {}={}&  b_{(1)}a_{(1)}\otimes  v _{(0)}\tl a_{(2)}\otimes \pi(b_{(2)}) v _{(1)}\pi(a_{(3)}).
    \label{eq:DeltaRHS}
\end{align}
Thus $A \otimes V$ is an object in~$\qMod{A}{A}{H}{A}$ and $A\otimes -$ defines a functor from~$\YD^H_A$ to~$\qMod{A}{A}{H}{A}$, which acts on morphisms in the obvious way.

Conversely, consider an $\mathcal{F} \in \qMod{A}{A}{H}{A}$ and recall the functor \mbox{$F:\qMod{A}{A}{}{A} \to \qMod{H}{}{}{A}$} introduced in \textsection \ref{section:TwoSidedHopfMods}. We  define  a right $H$-comodule structure on $F(\mathcal{F})$ by 
\begin{align*}
F(\mathcal{F}) \to F(\mathcal{F}) \otimes H, & & [f] \mapsto [f_{(2)}] \otimes \pi(f_{(1)})\pi(S(f_{(3)})),
\end{align*}
and hence give $F(\mathcal{F})$ a Yetter--Drinfeld structure. This gives a functor $F$ from  $\qMod{A}{A}{H}{A}$ to $\YD^H_A$, where morphisms are defined by descending to the quotient.

Finally, we note that an adjoint equivalence between the two categories is given by the unit natural transformation
\begin{align*}
\mathcal{F} \to A \otimes (F(\mathcal{F})), & & f \mapsto f_{(1)} \otimes [f_{(2)}],
\end{align*}
and the counit natural transformation 
\begin{align*}
  F (A \otimes \mathcal{F}) \to \mathcal{F}, & & a \otimes [f] \mapsto \e(a)f. 
\end{align*}

\begin{remark}\label{rem:YDHAbraiding}
In general, the formula for the Yetter\/--\/Drinfeld braiding~\eqref{eq:YDbraidingGeneral} is not
well defined for objects in~$\YD^H_A$. This difficulty can be avoided if~$\pi$ admits a bicovariant
splitting map~$i\colon H\to A$, as we will see below.
\end{remark}

\subsection{Left $A$-covariant right $H$-covariant first-order differential calculi}

If $(\Omega^1(A),\exd)$ is a left $A$-covariant right $H$-covariant first-order differential calculus
over~$A$ then $\Omega^1(A)$ is an object in~$\qMod{A}{A}{H}{A}$. In particular, since $(\Omega^1(A),\exd)$ is left $A$-covariant there is a corresponding ideal~$I_A$ of~$A^+$.
From~\cite[\textsection 2]{Maj} the right $H$-covariance condition is equivalent to the condition
$
  \Ad_\pi I_A \subseteq I_A\otimes H,
$
where
\begin{align}\label{eq:AdHpi}
  \Ad_\pi a := (\id\otimes\pi)\circ \Ad_R(a) = \ a_{(2)}\otimes \pi\big(S(a_{(1)})a_{(3)}\big),
  & & \text{for $a\in A$},
\end{align}
is the right $H$-coaction induced by the right adjoint $A$-coaction~\eqref{eq:AdA}.

\subsection{Restrictions of left $A$-covariant right $H$-covariant first-order differential calculi}\label{sec:ResDC}

Let $\Omega^1(A) \cong \Omega_u^1(A)/N_A$ be a left $A$-covariant right $H$-covariant first-order differential
calculus over~$A$ and $I_A = \xi_A(F(N_A))$ the corresponding $\Ad_{\pi}$-coinvariant ideal of~$A^+$.

For a quantum homogeneous space $B = A^{\co(H)}$, let $\Omega^1(B) :=  \mathrm{span}_{\mathbb{C}}\!\left\{ a\exd b\mid a,b\in B \right\}$ be the restriction of~$\Omega^1(A)$ to~$B$ and let
$I_B$ be the corresponding ideal of~$B^+$. Then  $I_B = I_A\cap B^+$, see~\cite[Theorem 5.77]{BeggsMajidBook} for details.

The adjoint right $H$-coaction $\Ad_\pi$, induced by~\eqref{eq:AdHpi}, acts on $b \in B$ as
\[
  \Ad_\pi b 
  = b_{(2)} \oby \pi(S(b_{(1)})b_{(3)}) = b_{(2)} \oby \pi(S(b_{(1)})),
\]
hence it coincides with the right $H$-coaction induced by the antipode from the left $H$-coaction
on~$B$. Therefore $I_B$ is a right $H$-subcomodule of~$I_A$. Note that we can define a right $H$-coaction
on~$\Phi(\Omega^1(B))$ using~$\Ad_\pi$ as follows
\begin{equation}\label{eq:YDcoact}
  \Ad_\pi [b] = [b_{(2)}] \otimes \pi(S(b_{(1)})b_{(3)}) = [b_{(2)}] \otimes \pi(S(b_{(1)})),\qquad\text{for $b\in B$}.
\end{equation}
We have  the following commutative diagram in the category~$\lMod{A}{}$
\begin{equation}\label{eq:diagiota}
  \begin{gathered} 
  \xymatrix{
    \Omega^1(A) \ar[rr]^{(\id\otimes\xi_A)\circ\unit} &\qquad & A\otimes F(\Omega^1(A)) \\
    \Omega^1(B) \ar[u]^{\iota^\prime}\ar[rr]^{(\id\otimes\xi_B)\circ\unit^\pi} &\qquad& A\square_H \Phi(\Omega^1(B))\ar[u]_{\id\otimes\iota}
  }
  \end{gathered}
\end{equation}
where $\iota$ and $\iota^\prime$ are the evident inclusions. It is important to note that $\iota$ is an $H$-comodule map.
In what follows we identity $\Phi(\Omega^1(B))$ with its image under~$\iota$.

In this paper, we will only deal with quantum principal bundles over quantum homogeneous spaces.
In this special case, we see that we have the following commutative diagram
\begin{align*}
\xymatrix{ 
\Omega_u^{1}(A)         \ar[d]_{\unit}           \ar[rrrd]^{\mathrm{ver}}       & & &      \\
  \,\, A \otimes A^+  \ar[rrr]_{\id \otimes \pi}                      & & & A \, \square_H H^+.  \\
}
\end{align*}
This means that for any subbimodule $N \sseq \Omega^1_u(A)$, with corresponding ideal $I := \unit(N)$, it holds that 
\[
\mathrm{can}(N) =  A \otimes \pi(I).
\]
From this we see that the requirement that $\mathrm{ver}(N) = A \oby I$ is automatically satisfied.

\begin{prop}[{\cite[\textsection 2]{Maj}}] Let $\pi\colon A\to H$ be a surjective Hopf algebra map such that $A$ is an
  $H$-Hopf--Galois extension of the associated quantum homogeneous space $B:=A^{\co(H)}$, and   $\Omega^1(A)$ is a left $A$-covariant right $H$-covariant differential caclulus over~$A$.
  Then $(A,\Delta_R,\Omega^1(A))$ is a quantum $H$-principal bundle.  
\end{prop}

\subsection{Constructing a Yetter--Drinfeld $H$-module structure}

In this section we establish a necessary list of conditions allowing us to associate a Yetter--Drinfeld module to a quantum $H$-principal bundle over a quantum homogeneous space. This is our principal tool for constructing a Nichols algebra presentation of the quantum Grassmannian Heckenberger--Kolb calculus in \textsection \ref{secDrinfeldJimbo}.

\begin{thm}\label{th:mainYD}
  Suppose $(A,\Delta_R,\Omega^1(A))$ is a quantum $H$-principal bundle over a quantum homogeneous space
  $B:=A^{\co(H)}$, and let $\Omega^1(B)$ be the restriction of~$\Omega^1(A)$ to~$B$.   Moreover, assume that
  
  (i) $\iota(\Phi(\Omega^1(B)))$ is a subobject of~$F(\Omega^1(A))$ in~$\Mod_A$,

  (ii) $\Omega^1(B)$ is a finite-dimensional object of~$\Modz{A}{B}$.
  
  \noindent Then 

  (a) a right $H$-action on~$\Phi(\Omega^1(B))$ is defined by 
  \begin{equation}\label{eq:YDact}
    \tl\colon \Phi(\Omega^1(B)) \oby H \to \Phi(\Omega^1(B)), \quad [b]\otimes h \mapsto [b]\tl_A i(h) = [b\, i(h)],
  \end{equation}
  where $i$ is a bicovariant splitting of~$\pi$, and $\tl_A$ is the right $A$-action on~$\iota(\Phi(\Omega^1(B)))$,

  (b) the triple $(\Phi(\Omega^1(B)),\tl,\Ad_\pi)$, where $\Ad_\pi$ is defined by~\eqref{eq:YDcoact},
  defines a Yetter--Drinfeld module over~$H$,
  
  (c) the Yetter--Drinfeld  module structure is independent of the choice of bicovariant splitting map~$i$.
\end{thm}

\begin{proof} For notational simplicity set $V := \Phi(\Omega^1(B))$.
For the $H$-action~\eqref{eq:YDact} to be well defined it should hold that 
\begin{equation}\label{action}
  v \tl_A\big(i(h h') - i(h)i(h')\big) = 0,\qquad\text{for all $v\in V$, $h,h'\in H$.}
\end{equation}
Now $\pi\big(i(h h') - i(h)i(h')\big) = 0$, and as discussed in \textsection \ref{sec:Tak}, it holds that $\ker \pi = B^+A$.
Therefore $i(h h') - i(h)i(h') \in B^+A$.  Hence, since $B^+$ acts trivially on~$V$, we must have that~\eqref{action} holds and
that $i$ defines an action of $H$ on $V$.

We now show that $\Ad_\pi$ and $\tl$ define a Yetter--Drinfeld structure on $V$. Note first that for $b\in B^+$
and $h\in H$ we have
\begin{align*}
 ([b]\tl h_{(2)})_{(0)} &{}\oby h_{(1)}([b]\tl h_{(2)})_{(1)}  \\
&{} = [b\, i(h_{(2)})]_{(0)} \oby h_{(1)}(b \, i(h_{(2)}))_{(1)}\\
&{} = [(b\, i(h_{(2)}) )_{(2)}] \oby
        h_{(1)}\pi\Big(S\big( (b\, i(h_{(2)}))_{(1)}\big) (b\, i(h_{(2)}))_{(3)}\Big)\\
  &{} = [b_{(2)} i(h_{(2)})_{(2)}]
    \otimes h_{(1)} S\Big(\pi\big(i(h_{(2)})_{(1)}\big)\Big)\pi\big(S(b_{(1)})b_{(3)}\big)\pi\big(i(h_{(2)})_{(3)}\big).
\end{align*}
Since $i$ is a bicovariant splitting map we have that
\[
  \pi(i(h)_{(1)})\otimes i(h)_{(2)} \otimes \pi(i(h)_{(3)}) = h_{(1)}\otimes i(h_{(2)})\otimes h_{(3)}.
\]
Combining this with the previous calculation we see that
\begin{align*}
  ([b]\tl h_{(2)})_{(0)} \oby h_{(1)}([b]\tl h_{(2)})_{(1)}
  & {} = [b_{(2)}i(h_{(3)})] \oby h_{(1)} S(h_{(2)})\pi\big(S(b_{(1)})b_{(3)}\big)h_{(4)}\\
  & {} = [b_{(2)}] \tl h_{(1)} \oby \pi\big(S(b_{(1)})b_{(3)}\big) h_{(2)}\\
  & {} = [b]_{(0)} \tl h_{(1)} \oby [b]_{(1)} h_{(2)},
\end{align*}
as required. It follows that the action and coaction satisfy the Yetter--Drinfeld condition~\eqref{eq:YDcondition}.

Finally, we show that the Yetter--Drinfeld structure is independent of the choice of bicovariant splitting map. Let $i'$ be a second splitting map and note that since
\[
  \pi(i(h) - i'(h)) = \pi \circ i(h) - \pi \circ i'(h) = h - h = 0,
\]
we must have
\(
\mathrm{Im}(i-i') \sseq \ker(\pi) = B^+A
\). In particular, $i(h) - i'(h) \in B^+A$ and $\varepsilon( i(h)-i'(h))=0$.
Since we have assumed that~$\Omega^1(B)$ is an object in~$\Modz{H}{B}$,  the right $B$-action on~$V$ is trivial.
Therefore
\[
  [b] \tl_A i(h) - [b] \tl_A i'(h)
  = [b] \tl_A (i(h) - i'(h))
  = \varepsilon (i(h) - i'(h)) [b]
  = 0,
\]
implying that the two actions are equal.
\end{proof}

\begin{remark}
Recall from \S\ref{sec:Tak} that $\Psi$ is a monoidal functor. Thus the morphism $\Psi(\sigma)$
satisfies the Yang--Baxter equation for linear operators on  $\Omega^1(B)\otimes_B\Omega^1(B)$.
Thus we can consider the corresponding Nichols algebra~$\fB(\Omega^1)$ as a quotient of the tensor
algebra~$\cT_B(\Omega^1(B))$ of~$\Omega^1(B)$ over~$B$, see~\cite{WoronowiczDC89}.  This means that $\fB(\Omega^1(B)) \cong \Psi(\fB(V))$, giving us a Nichols algebra description of the  differential calculus.
\end{remark}

\begin{remark}\label{rem:NicholsAH} If $i\colon H \to A$ is a bicovariant splitting then
  $\qMod{A}{A}{H}{A}\cong\YD^H_A$ and  so the corresponding Nichols algebra makes sense with respect to the braiding
  \begin{equation}\label{eq:braidingAH}
    \sigma(v\otimes w) = w_{(0)} \otimes v (\tl_A\circ i)(w_{(1)}).
  \end{equation}
  
\end{remark}

\section{The quantum Grassmannian Heckenberger--Kolb calculi}
\label{sec:Grass}

In this section we apply Theorem~\ref{th:mainYD} to the first-order parts of the Heckenberger--Kolb calculi
of the quantum Grassmannians to produce a Yetter--Drinfeld structure on their cotangent spaces. We then build on this construction to present the holomorphic and anti-holomorphic subcomplexes as Nichols algebras.

\subsection{Preliminaries on Drinfeld--Jimbo quantum groups} \label{secDrinfeldJimbo}

In this section we recall basic material about Drinfeld--Jimbo quantised universal enveloping algebras~\cite{DrinfeldICM,Jimbo1986} and
their representation theory. We refer the reader to~\cite{KSLeabh,ChariPressley,MajLeabh} for further details.

\subsubsection{Drinfeld--Jimbo quantised universal enveloping algebras}\label{sec:DJ}
Let $\frak{g}$ be a finite-dimensional complex semi\-simple Lie algebra of rank~$r$. We fix a~Cartan subalgebra~$\frak{h}$ with corresponding root system $\Delta \sseq \frak{h}^*$, where~$\frak{h}^*$ denotes the linear dual of~$\frak{h}$.  
Fix a choice of simple roots $\{\alpha_1, \dots, \alpha_r\}$.
Denote by $(\cdot,\cdot)$ the symmetric bilinear form induced on $\frak{h}^*$ by the  Killing form of
$\frak{g}$, normalised so that any shortest simple root $\alpha_i$ satisfies $(\alpha_i,\alpha_i) = 2$.
Let $\{\varpi_1, \dots, \varpi_r\}$ denote the corresponding set of fundamental weights of~$\mathfrak{g}$. 
The Cartan matrix $A = (a_{ij})$ of $\frak{g}$ is the $(r \times r)$-matrix defined by
\(
a_{ij} := \big(\alpha_i^{\vee},\alpha_j\big),
\)
where $\alpha_i^{\vee} := 2\alpha_i/(\alpha_i,\alpha_i)$. 
Let  $q \in \bC$ such that $q$ is not a~root of unity and denote $q_i := q^{(\alpha_i, \alpha_i)/2}$. The \emph{quantised enveloping algebra}  $U_q(\frak{g})$ is the  noncommutative associative  algebra  generated by the elements   $E_i, F_i, K_i$, and  $K^{-1}_i$, for $ i=1, \ldots, r$,  subject to the relations 
\begin{gather*}
 K_iE_j =  q_i^{a_{ij}} E_j K_i, \quad  K_iF_j= q_i^{-a_{ij}} F_j K_i, \quad  K_i K_j = K_j K_i, \quad K_iK_i^{-1} = K_i^{-1}K_i = 1,\\
  E_iF_j - F_jE_i  = \delta_{ij}\frac{K_i - K\inv_{i}}{q_i-q_i^{-1}},
\end{gather*}
along with the \emph{quantum Serre relations}  
\begin{align*}
  \sum\nolimits_{s=0}^{1-a_{ij}} (-1)^s  \begin{bmatrix} 1 - a_{ij} \\ s \end{bmatrix}_{q_i}
   E_i^{1-a_{ij}-s} E_j E_i^s = 0,\quad \textrm{ for }  i\neq j,\\
  \sum\nolimits_{s=0}^{1-a_{ij}} (-1)^s \begin{bmatrix} 1 - a_{ij} \\ s \end{bmatrix}_{q_i}
   F_i^{1-a_{ij}-s} F_j F_i^s = 0,\quad \textrm{ for }  i\neq j,
\end{align*}
where we have used the $q$-binomial coefficients defined as follows
\begin{gather*}
  [n]_q! := [n]_q[n-1]_q\ldots [2]_q[1]_q,\qquad
  \text{where $[k]_q := \frac{q^k-q^{-k}}{q-q^{-1}}$,}\\
  \begin{bmatrix} n \\ k \end{bmatrix}_{q} := \frac{[n]_q!}{[k]_q![n-k]_q!}.
\end{gather*}
A Hopf algebra structure is defined  on $U_q(\frak{g})$ by
\begin{align*}
\Delta(K_i) &= K_i \oby K_i,\quad  \Delta(E_i) = E_i \oby K_i + 1 \oby E_i, \quad \Delta(F_i) = F_i \oby 1 + K_i\inv  \oby F_i,\\
& \qquad S(E_i) =  - E_iK_i\inv,    \quad S(F_i) =  -K_iF_i, \quad  S(K_i) = K_i\inv,  \\
&\qquad \qquad \qquad \e(E_i) = \e(F_i) = 0, ~~ \e(K_i) = 1.     
\end{align*}  
Let $\mathcal{P}$ be the weight lattice of~$\fg$, and $\mathcal{P}^+$ its set of dominant integral weights.
We consider~$\Rep_1U_q(\fg)$,
the full subcategory of the category of (left) $U_q(\fg)$-modules,
whose the objects are finite-dimensional $U_q(\fg)$-modules having a~weight
decomposition
\(
V = \bigoplus_{\mu\in\mathcal{P}}V(\mu)
\).
Recall that a vector $v\in V$ is called a \emph{weight vector} of \emph{weight}~$\mu\in\fh^*$ if
$K_i\tr v = q^{(\alpha_i,\mu)}v$ for all $i=1,\ldots,r$.
The category~$\Rep_1U_q(\fg)$ is a~semisimple tensor category whose simple objects are irreducible modules~$V_\lambda$ with
 highest weight~$\lambda\in\mathcal{P}^+$.
The character of~$V_\lambda$ is given by the
classical Weyl character formula for the irreducible $\fg$-module~$\what V_\lambda$ with highest
weight~$\lambda$. In fact,
the category~$\Rep_1U_q(\fg)$ is equivalent to the category~$\cO_f$ of finite-dimensional representations
of~$\fg$.
We refer to~\cite[\S5.8]{EtingofTensorCat} and~\cite[\S7]{KSLeabh} for further details.

\subsection{Quantum coordinate algebras} 
In this subsection we recall some necessary material about quantised coordinate algebras, see~\cite[\textsection6 and \textsection7]{KSLeabh}
and~\cite{FRT} for further details.
Let $V$ be a finite-dimensional left $U_q(\frak{g})$-module, $v \in V$, and $f \in V^*$, where $V^*$ is the $\mathbb{C}$-linear dual of $V$ endowed with its  right \mbox{$U_q(\frak{g})$-module} structure. An~important point to note is that, with respect to the equivalence of left and right $U_q(\frak{g})$-modules given by the invertible antipode, the left module corresponding to $V^*_{\mu}$ is isomorphic to $V_{-w_0(\mu)}$, where $w_0$ denotes the longest element in the Weyl group of $\frak{g}$.

Consider the function  $c^{V}_{f,v}:U_q(\frak{g}) \to \bC$ defined by $c^{V}_{f,v}(X) := f\big(X \triangleright v\big)$. The {\em coordinate ring} of $V$ is the subspace
\[
C(V) := \textrm{span}_{\mathbb{C}}\!\left\{ c^{V}_{f,v} \,| \, v \in V, \, f \in V^*\right\} \sseq U_q(\frak{g})^\circ,
\]
where $U_q(\frak{g})^\circ$ denotes the Hopf dual of $U_q(\frak{g})$.
A $U_q(\fg)$-bimodule structure on~$C(V)$ is given by
\begin{equation} \label{eq:Zact}
  (Y\triangleright c^{V}_{f,v} \triangleleft Z)(X) := f\left((ZXY)\triangleright v\right)
  = c^{V}_{f\triangleleft Z, Y \triangleright v}  (X).
\end{equation}
 It is easily checked that $C(V) \sseq U_q(\frak{g})^\circ$, and moreover that a Hopf subalgebra of $U_q(\frak{g})^\circ$ is given by 
\begin{equation}\label{eq:PeterWeyl}
  \cO_q(G) := \bigoplus_{\mu \in \mathcal{P}^+} C(V_{\mu}).
\end{equation}
We call $\cO_q(G)$ the {\em quantum coordinate algebra of~$G$}, where~$G$ is the compact, connected, simply-connected, simple Lie group  having~$\frak{g}$ as its complexified Lie algebra.

Let $H$ be a Hopf algebra and $V$ a right $H$-comodule. Suppose that $U$ is a Hopf algebra which is dualy
paired with~$H$ via a pairing $\langle\cdot,\cdot\rangle\colon U\otimes H \to \Cee$. Recall that a right
$U$-action on~$V$ is defined by
\begin{equation}\label{eq:coactTOact}
  \ltr\colon  V \otimes U \to V,\qquad
  v\otimes X \mapsto v_{(0)}\langle S(X), v_{(1)} \rangle.
\end{equation}
Similarly one may define a left action of~$U$ on a left $H$-comodule. 

\subsubsection{Quantum exterior algebras}\label{sec:qExtAlg}

The category~$\Rep_1 U_q(\fg)$ is a braided monoidal category where the braiding comes
from the universal $R$-matrix~$R\in U_q(\fg)\widehat{\otimes}U_q(\fg)$ as defined in~\cite{DrinfeldICM}. For any two objects $V$, $W \in\Rep_1 U_q(\fg)$,
let 
\begin{align*}
\rho_V\colon U_q(\fg)\to \End(V), & & \rho_W\colon U_q(\fg)\to \End(W),
\end{align*}
be the corresponding structure maps.
The braiding~$\sigma_{R,V\otimes W}\colon V\otimes W\to W\otimes V$ induced by~$R$
is defined by
\begin{equation}\label{eq:RmatBraiding}
  \sigma_{R,V\otimes W} := \tau \circ (\rho_V\otimes\rho_W)(R),
\end{equation}
where $\tau\colon V\otimes W \to W\otimes V$ is the ordinary flip. Following~\cite{DrinfeldQuasiHopf}, define the \emph{normalised braiding} $\tilde{\sigma}_{R,V\otimes W}\colon
V\otimes W\to W\otimes V$ by
\begin{equation}\label{eq:RmatBraidingNorm}
  \tilde{\sigma}_{R,V\otimes W} := \sqrt{\sigma_{R,W\otimes V}^{-1}\sigma_{R,V\otimes W}^{-1}}\, \sigma_{R,V\otimes W}.
\end{equation}
\begin{remark} 
   Note that formula~\eqref{eq:RmatBraidingNorm} is not well defined for an arbitrary braided
  monoidal category. However, it is well defined for the quantised universal enveloping algebra~$U_q(\fg)$ since
  the universal $R$-matrix $R\in U_q(\fg)\what{\otimes} U_q(\fg)$ can be decomposed as
  \[
    R = R_0R_1 = R_1R_0,
  \]
  where $R_0$ is ``the diagonal part''  of~$R$, and $R_1$ is unipotent. We refer
  to~\cite[\S3]{DrinfeldQuasiHopf} for further discussion.
\end{remark}

As in~\S\ref{sec:braidedHopf} we write~$\tilde{\sigma}_{R}$ instead of~$\tilde{\sigma}_{R,V\otimes W}$.
Since  $\tilde{\sigma}_R^2  =\id$, $\tilde{\sigma}_R$ is a~symmetric commutativity constraint but it does not satisfy  the Yang\/--Baxter equation in general, see~\cite[\S3]{DrinfeldQuasiHopf} and~\cite{BerensteinZwicknagl2008}.

For any $V\in\Rep_1(U_q(\fg))$, denote
\begin{equation}\label{eq:qSLambda}
  S^2_qV := \{ x\in V\otimes V\mid \tilde{\sigma}_R(x) = x\},\quad
  \Lambda^2_qV := \{ x\in V\otimes V\mid \tilde{\sigma}_R(x) = -x\}.
\end{equation}
Note that $V\otimes V = S^2_q V \oplus \Lambda^2_q V$.
Following~\cite{BerensteinZwicknagl2008},  define
the \emph{quantum exterior algebra} $\Lambda_q(V)$ of~$V$ to be
\[
  \Lambda_q(V) := \mathcal{T}(V) / \langle S^2_qV \rangle,
\]
where $\mathcal{T}(V)$ is the tensor algebra of~$V$ and $\langle I \rangle$ denotes the two-sided ideal
generated by $I\subset\mathcal{T}(V)$.
A $U_q(\fg)$-module~$V$ is called \emph{flat} if the Hilbert\/--\/Poincar\'e series of~$\Lambda_q(V)$ is the
same as for its classical counterpart~$\Lambda(\what{V})$.

\begin{example}\label{ex:VV}
Fix two positive integers $n$ and $m$ such that $n>m$.
Let $V$ be a~$U_q(\fgl_m)$-module and $W$ a~$U_q(\fsl_{n-m})$-module. Denote by $V\boxtimes W$ the
tensor product~$V\otimes W$ viewed as a~$U_q(\fgl_m\oplus\fsl_{n-m})$-module.
We have the following isomorphism of~$U_q(\fgl_m\oplus\fsl_{n-m})$-modules
\[
  (V_{\varpi_1}\boxtimes V_{\varpi_1}^*)\otimes (V_{\varpi_1}\boxtimes V_{\varpi_1}^*)
  \simeq  (V_{\varpi_2}\boxtimes V_{2\varpi_1}^*) \oplus
  (V_{\varpi_2}\boxtimes V_{\varpi_2}^*) \oplus
  (V_{2\varpi_1}\boxtimes V_{2\varpi_1}^*) \oplus
  (V_{2\varpi_1}\boxtimes V_{\varpi_2}^*).
\]
For details see Table~5 in~\cite{VinbergOnishchik}.

As was shown in~\cite[Proposition~2.33]{BerensteinZwicknagl2008}, we have
\[
  S^2_q(V_{\varpi_1}\boxtimes V_{\varpi_1}^*) =   (V_{\varpi_2}\boxtimes V_{\varpi_2}^*) \oplus
  (V_{2\varpi_1}\boxtimes V_{2\varpi_1}^*).
\]
Moreover, the quantum exterior
algebra~$\Lambda_q(V_{\varpi_1}\boxtimes V_{\varpi_1}^*)$ has the same Hilbert\/--Poincar\'e series as the exterior
algebra~$\Lambda(\what V_{\varpi_1}\boxtimes \what V_{\varpi_1}^*)$.
Hence $V_{\varpi_1}\boxtimes V_{\varpi_1}^*$ is flat.
\end{example}

\subsubsection{Quantum Howe duality}\label{sec:qHowe}

In~\cite[Theorem~4.2.2]{CautisKamnitzerMorrison2014} (see also~\cite{Zhang2002}),
it was shown that
as $U_q(\fgl_m \oplus \fsl_{n-m})$-modules the quantum exterior algebra of $V_{\varpi_1}\boxtimes V_{\varpi_1}^*$ decomposes
as 
\begin{equation}\label{eq:qHowe}
  \Lambda_q( V_{\varpi_1} \boxtimes V_{\varpi_1}^* ) = \bigoplus_{\lambda} V_\lambda \boxtimes V^*_{\lambda^t},
\end{equation}
where $\lambda$ varies over all $(n-m)$-bounded partitions, and $\lambda^t$ is the
transpose of~$\lambda$.

\subsection{Preliminaries on quantum flag manifolds}\label{sec:qflags}

Let $\{\alpha_i\}_{i\in S}$ be a~subset of simple roots. In what follows, by abuse of notation, we denote
by~$S$ not only an~index subset but also the corresponding subset of simple roots~$\{\alpha_i\}_{i\in S}$.
Consider the Hopf 
subalgebra
\[
  U_q(\frak{l}_S) := \left< K_i, E_j, F_j \mid i = 1, \ldots, r; j \in S \right>.
\]
The Hopf 
algebra embedding $\iota_S:U_q(\frak{l}_S) \hookrightarrow U_q(\frak{g})$ induces a~dual Hopf 
algebra map $\iota_S^{\circ}: U_q(\frak{g})^{\circ} \to U_q(\frak{l}_S)^{\circ}$.
By construction $\cO_q(G) \sseq U_q(\frak{g})^{\circ}$, so the restriction map
\begin{align}\label{eq:piS}
\pi_S:= \iota_S^{\circ}|_{\cO_q(G)}: \cO_q(G) \to U_q(\frak{l}_S)^{\circ},\
\end{align}
defines a Hopf 
subalgebra
\(
\cO_q(L_S) := \pi_S\big(\cO_q(G)\big) \sseq U_q(\frak{l}_S)^\circ.
\)
The {\em quantum flag manifold} associated to $S$ is the quantum homogeneous space associated to the
surjective Hopf 
algebra map  $\pi_S:\cO_q(G) \to \cO_q(L_S)$,  and is denoted by
\[
  \cO_q\big(G/L_S\big) := \cO_q \big(G\big)^{\co\left(\cO_q(L_S)\right)}.
\]

\begin{remark}\label{rem:splittingmap}
  Since $\cO_q(L_S)$ is cosemisimple, it follows from Lemma~\ref{cosemiH} that the extension
  $\cO_q(G/L_S) \hookrightarrow \cO_q(G)$ is a principal comodule algebra and~$\pi_S$ admits a~bicovariant
  splitting map.
\end{remark}

\subsubsection{The Heckenberger--Kolb calculus}

Let $S$ be a subset a simple roots of~$\fg$.  Following the classical case (see, for example,~\cite{BastonEastwood})
we say that the quantum flag manifold associated to~$S$ is of \emph{irreducible type} if
$\fg/\fl_S$ is a~direct sum of two dual irreducible $\fl_S$-modules. The following theorem summarises~\cite[Theorem~7.2]{HK} and~\cite[Propositions~3.6 and~3.7]{HKdR}.

\begin{thm}\label{th:HK}
  For any irreducible quantum flag manifold $\cO_q(G/L_S)$, there exist
  exactly two non-isomorphic, irreducible, left $\cO_q(G)$-covariant, finite-dimensional, first-order
  differential calculi
  \[
    \Omega_q^{(1,0)}(G/L_S),\ \Omega_q^{(0,1)}(G/L_S) \in \Modz{\cO_q(G)}{\cO_q(G/L_S)},
  \]
  and the corresponding maximal prolongations
  $\Omega_q^{(\bullet,0)}(G/L_S)$ and~$\Omega_q^{(0,\bullet)}(G/L_S)$
  have classical dimension.
\end{thm}
\begin{remark}
  For the reader's convenience, we recall some detail about the proof of Theorem~\ref{th:HK}. The
  left-covariant first-order differential calculi were first classified in \cite{HK}, where it was shown that,
  up to isomorphism, there exist precisely two irreducible finite-dimensional left-covariant differential
  calculi over the irreducible quantum flag manifolds. This was achieved by classifying the equivalent
  notation of a quantum tangent space using the coradical filtration of the locally finite part of the dual
  coalgebra of $\cO_q(G/L_S)$. The maximal prolongation of these first-order calculi to differential calculi
  was then explicitly described in \cite{HKdR} where it was shown, among many other things, that the calculi
  have classical dimension.
\end{remark}

\subsection{Quantum Grassmannians}
Consider the $A$-series irreducible quantum flag manifolds, namely, the quantum
Grassmannians. Let $\fg = \fsl_{n}$ and $m$ be an~integer such that $1\leq m < n$. Fix
$S = \{ \alpha_i \mid i\in \{1,\ldots,n-1\}\setminus\{m\}\}$ a~subset of simple roots of~$\fg$. In this case
$\mathfrak{l}_S = \fgl_m\oplus\fsl_{n-m}$. The quantum flag manifold associated to~$S$ is called \emph{the
  quantum $(n,m)$-Grassmannian} and denoted by~$\cO_q(\rGr_{n,m})$. We denote by
$\pi := \pi_S \colon \cO_q(SU_{n})\to \cO_q(U_m\times SU_{n-m})$
the Hopf algebra map corresponding to~$S$~\eqref{eq:piS}.

In what follows we use the fact that the quantum coordinate algebra~$\cO_q(SU_{n})$ can be described in terms
of matrix coefficients of the first fundamental representation~$V_{\varpi_1}$ of~$U_q(\fsl_{n})$, recalling the
Faddeev\/--Reshetikhin\/--Takhtadzhyan approach~\cite{FRT}. In particular, the matrix coefficients 
$u^i_j:=c^{\varpi_1}_{f_i,v_j}$, $i,j=1,\ldots,n$ generate~$\cO_q(SU_{n})$.
The corresponding dual pairing between~$U_q(\fsl_{n})$ and $\cO_q(SU_{n})$ is given by
\begin{equation}\label{eqn:dualpairing}
  \langle K_i,u^j_j\rangle = q^{\delta_{i,j-1} - \delta_{i,j}},\qquad
  \langle E_i,u^{i+1}_i\rangle = 1,\qquad
  \langle F_i,u^{i}_{i+1}\rangle = 1.
\end{equation}
Therefore we can view any
$\cO_q(SU_{n})$-comodule as a~$U_q(\fsl_{n})$-module via~\eqref{eq:coactTOact}.

\subsection{Quantum principal bundles over the quantum Grassmannians}

In this section we recall necessary facts about the first-order part of the Heckenberger--Kolb calculi over
the quantum Grassmannians. These allow us to apply Theorem~\ref{th:mainYD} to construct the associated Yetter--Drinfeld modules. Throughout we denote by  $\Omega^1_q(\rGr_{n,m}) := \Omega_q^{(1,0)}(\rGr_{n,m}) \oplus \Omega_q^{(0,1)}(\rGr_{n,m})$ the
direct sum of the pair of Heckenberger--Kolb calculi over~$\cO_q(\rGr_{n,m})$.

\begin{prop}[{\cite[Corollary~5.3 and Proposition~5.8]{BWGrass}}]
  \label{prop:GrCaluli}
Let $\Omega^1_{bc}(SU_{n})$ be the bicovariant first-order differential calculus over~$\cO_q(SU_{n})$
  associated to the $R$-matrix of~$U_q(\fsl_{n})$.  
 For every $m=1,\ldots,n-1$, there is an ideal~$N_m$ of~$\Omega_{bc}^1(SU_{n})$ such that
  
  1) the quotient calculus $\Omega^1_m(SU_{n}):=\Omega_{bc}^1(SU_{n})/N_m$ is a left $\cO_q(SU_{n})$-covariant right
  $\cO_q(U_m\times SU_{n-m})$-covariant first order differential calculus,
  
  2) the restriction of~$\Omega_m^1(SU_{n})$ to~$\cO_q(\mathrm{Gr}_{n,m})$ is~$\Omega^1_q(\rGr_{n,m})$.
\end{prop}

Set $V^{(1,0)} := \Phi(\Omega_q^{(1,0)}(\rGr_{n,m}))$ and~$V^{(0,1)} := \Phi(\Omega_q^{(0,1)}(\rGr_{n,m}))$.

\begin{cor}\label{cor:VisYD}
  The quantum principal bundle above determines $\cO_q(U_m\times SU_{n-m})$-Yetter--Drinfeld module
  structures for~$V^{(1,0)}$ and~$V^{(0,1)}$.
\end{cor}
\begin{proof}
  As was noted in Remark~\ref{rem:splittingmap}, $\pi$ admits a bicovariant splitting map.
  From Proposition~\ref{prop:GrCaluli} it follows that the triple
  $(\cO_q(SU_{n}), \cO_q^1(\mathrm{Gr}_{n,m}),  \Omega^1_m(SU_{n}))$
  defines a quantum principal bundle.
  By Theorem~\ref{th:HK}, we have that $\Omega_q^{(1,0)}(\mathrm{Gr}_{n,m})$ and $\Omega_q^{(0,1)}(\mathrm{Gr}_{n,m})$
  are in~$\Modz{\cO_q(G)}{\cO_q(G/L_S)}$.
  Moreover, $\iota(V^{(1,0)})$ and~$\iota(V^{(0,1)})$ are right $\cO_q(SU_{n})$-submodules,
  see~\eqref{eq:V01YDact} below.
  Thus the conditions of Theorem~\ref{th:mainYD} are satisfied.
\end{proof}

As noted in Remark~\ref{rem:NicholsAH}, since $\Omega^1_m(SU_n)$ is a left $\cO_q(SU_n)$-covariant and 
right $\cO_q(U_m\times SU_{n-m})$-covariant first order differential calculus, it admits a~braiding
by~\eqref{eq:braidingAH}. Determining the properties of the corresponding Nichols algebras presents itself as an interesting direct of research.

\subsection{The Yetter--Drinfeld braiding on $V^{(0,1)} \otimes V^{(0,1)}$}

Recall from \S\ref{sec:ResDC} that we have an embedding of right $H$-comodules
\[
  \iota \colon V^{(1,0)}\oplus V^{(0,1)} \to F(\Omega^1_m(SU_{n})).
\]
Let
\(
M := \{1, \dots, m\}
\), and
\(
  \barM := \{m+1, \dots, n\}
\).
From~\cite[Lemma~5.7]{BWGrass} it holds that
\begin{enumerate}
\item the set $\{[u^i_j] \mid  (i,j) \in \barM \by M\}$ is a basis of~$\iota(V^{(1,0)})$,
\item the set  $\{[u^i_j] \mid (i,j) \in M \by \barM\}$ is a basis of~$\iota(V^{(0,1)})$.
\end{enumerate}
In what follows, we will use the bases of~$V^{(1,0)}$ and $V^{(0,1)}$ induced by the bases
of~$\iota(V^{(1,0)})$ and~$\iota(V^{(0,1)})$ from the previous proposition.
The $\cO_q(U_m\times SU_{n-m})$-Yetter--Drinfeld module structure on $V^{(0,1)}$  is explicitly described below in terms of these bases.

As was shown in~\cite[Proposition~B.3]{BWGrass},
the $\cO_q(U_m\times SU_{n-m})$-action on~$V^{(1,0)}$
and $V^{(0,1)}$ is given as follows.  
  For $i\neq j$ and $(p,s)\in (M\times M) \cup (\barM\times \barM)$ we have
\begin{align}
\label{eq:V01YDact}
   &[u^i_j]\tl \pi(u^p_p) =  q^{-2/n}q^{\delta_{p,i}+\delta_{p,j}}[u^i_j],\\
   &[u^i_j]\tl \pi(u^p_s) =  q^{-2/n}\left(\theta(p-s)\delta_{s,i}(q-q^{-1})[u^p_j]
     + \theta(s-p)\delta_{p,j}(q-q^{-1})[u^i_s]\right)
    \nonumber
      \\
   &[u^i_j]\tl \pi(S(u^p_p)) = q^{2/n}q^{-\delta_{p,i}-\delta_{p,j}}[u^i_j],
    \nonumber
      \\
   &[u^i_j]\tl \pi(S(u^p_s)) =
      q^{2/n}\left( \delta_{s,i}\theta(p-s)(q^{-1}-q)[u^p_j] + \delta_{p,j}\theta(s-p)q^{2(p-s)}(q^{-1}-q)[u^i_s]
       \right). \nonumber
     \end{align}

 By~\eqref{eq:YDcoact}, the right $\cO_q(U_m\times SU_{n-m})$-coactions on~$V^{(1,0)}$ and~$V^{(0,1)}$ are
 given by
 \begin{equation}\label{eq:V01YDcoact}
   \begin{aligned} 
     \Ad_\pi[u^i_j]& = \sum_{(a,b)\in M \by \barM} [u^a_b] \oby \pi(S(u^i_a)u^b_j), & \text{ for } (i,j) \in M \by \barM\\
     \Ad_\pi[u^i_j] &= \sum_{(a,b) \in \barM \by M} [u^a_b] \oby \pi(u^b_j S(u^i_a)), & \text{ for } (i,j) \in \barM \by M. 
   \end{aligned}
\end{equation}

 \begin{lemma}\label{lem:YDbraiding01}
For  $[u^i_j] \otimes [u^k_l] \in V^{(0,1)} \oby V^{(0,1)}$, and $\sigma$ the Yetter--Drinfeld braiding
of~$V\ahol$, it holds that 
\begin{align}\label{eq:YDbraiding}
\sigma([u^i_j]&{} \otimes [u^k_l]) = q^{-\delta_{ki} + \delta_{lj}} [u^k_l] \otimes [u^i_j] + q^{-\delta_{ik}}(q - q^{-1}) \theta(l-j) [u^k_j] \otimes [u^i_l] \notag\\
&-  q^{\delta_{lj}} (q - q^{-1}) \theta(k-i)[u^i_l ] \otimes [u^k_j] 
- (q-q^{-1})^2 \theta(l-j)  \theta(k-i) [u^i_j] \otimes [u^k_l].
\end{align}
\end{lemma}

\begin{proof}
  Applying~\eqref{eq:V01YDcoact}  to~\eqref{eq:YDbraidingGeneral}, 
  we see that the non-zero terms are 
\begin{align*}
  \sigma\big([u^i_j]
  {} \oby [u^k_l]\big) ={}
  &  [u^k_l] \oby \left([u^i_j]\tl \pi(S(u^k_k)u^l_l)\right) + [u^k_j] \oby\left( [u^i_j]\tl\pi(S(u^k_k)u^j_l)\right)\\
  &{}  +  [u^i_l] \oby \left([u^i_j]\tl\pi(S(u^k_i)u^l_l)\right) +  [u^i_j] \oby \left([u^i_j]\tl\pi(S(u^k_i)u^j_l)\right).  
 \end{align*} 
A further application of~\eqref{eq:V01YDact} gives us~\eqref{eq:YDbraiding}.
\end{proof}

\subsection{The spectrum of the Yetter--Drinfeld braiding on $V\ahol \otimes V\ahol$}
Following~\eqref{eq:coactTOact}, we convert the $\cO_q(U_m\times SU_{n-m})$-coaction on~$V^{(0,1)}$ to a~$U_q(\fgl_m\oplus\fsl_{n-m})$-action.
A direct computation confirms that  the highest weight vectors in the $U_q(\fgl_m\oplus\fsl_{n-m})$-module $V\ahol\oby V\ahol$ are given by
  \begin{align} \label{HWVinV01}
    v_1 :={}&{} [u^1_n]\otimes [u^1_n],\\
    v_2 :={}&{} -q [u^1_n]\otimes [u^1_{n-1}] + [u^1_{n-1}]\otimes [u^1_n], \nonumber \\
    v_3 :={}&{} -q [u^1_n]\otimes [u^2_n] + [u^2_n]\otimes [u^1_n], \nonumber \\
    v_4 :={}&{} q^2 [u^1_n]\otimes [u^2_{n-1}] + [u^2_{n-1}]\otimes [u^1_n] - q ([u^1_{n-1}]\otimes [u^2_n] + [u^2_n]\otimes [u^1_{n-1}]). \nonumber
  \end{align}
Recall from~\cite{HK} that $V^{(0,1)} \simeq V_{\varpi_1}\boxtimes V_{\varpi_1}^*$. In terms of the
decomposition given in Example~\ref{ex:VV}, we have that
$v_1$ is the highest weight vector of~$V_{\varpi_2}\boxtimes V_{\varpi_2}^*$,
that $v_2$ is the highest weight vector of~$V_{\varpi_2}\boxtimes V_{2\varpi_1}^*$,
that $v_3$ is the highest weight vector of~$V_{2\varpi_1}\boxtimes V_{\varpi_2}^*$,
and that $v_4$ is the highest weight vector of~$V_{2\varpi_1}\boxtimes V_{2\varpi_1}^*$.

\begin{lemma}\label{lemma:YDbraidingOnHWV01}
  The braiding~\eqref{eq:YDbraiding} acts on the highest weight vectors~\eqref{HWVinV01}
  as follows:
  \[
    \sigma(v_1) = v_1,\quad \sigma(v_2) = -q^{-2} v_2,\quad \sigma(v_3) = - q^2 v_3,\quad \sigma(v_4)=v_4.
  \]
\end{lemma}

 Since we are interested in the $q$-deformed exterior algebra, it is natural to consider the Nichols algebra
given by the rescaled braiding~$-\sigma$; see~\cite[\S3]{WoronowiczDC89}.

\begin{cor}\label{cor:kerS2}
  It holds that
  \[
    \ker(\mathfrak{S}_2^{-\sigma}) =
    V_{\varpi_2}\boxtimes V_{\varpi_2}^* \oplus V_{2\varpi_1}\boxtimes V_{2\varpi_1}^*
    = S_q^2(V^{(0,1)}).
  \]
\end{cor}

\subsection{A Nichols algebra presentation of~$\Omega_q^{(\bullet,0)}(\rGr_{n,m})$ and~$\Omega_q^{(0,\bullet)}(\rGr_{n,m})$}
In this section we present the main result of the paper. Namely, we show that the holomorphic and
anti-holomorphic quantum exterior algebras of the subcomplex $\Omega_q^{(\bullet,0)}(\rGr_{n,m})$
and the subcomplex $\Omega_q^{(0,\bullet)}(\rGr_{n,m})$ are Nichols algebras.

\begin{thm}\label{th:mainNichols}
  There exist $\cO_q(U_m\times SU_{n-m})$-comodule algebra isomorphisms,
  or equivalently, $U_q(\fgl_{m}\oplus\fsl_{n-m})$-module algebra isomorphisms,
  \[
    \fB(V^{(1,0)},-\sigma) \simeq \Phi(\Omega_q^{(\bullet,0)}(\rGr_{n,m})),\qquad
    \fB(V^{(0,1)},-\sigma) \simeq \Phi(\Omega_q^{(0,\bullet)}(\rGr_{n,m})),
  \]
  where $\sigma$ is the Yetter--Drinfeld braiding~\eqref{eq:YDbraiding}.
\end{thm}
\begin{proof}
  Let $V := \Phi(\Omega^{(0,1)}(\rGr_{n,m}))$ and $N := \dim V = m(n-m)$.
  Recall from~\cite[\S6]{HK} that, as $U_q(\fgl_m\oplus\fsl_{n-m})$-modules, 
  \[
    V \simeq V_{\varpi_1}\boxtimes V_{\varpi_1}^*\qquad\text{and}\qquad
    \Phi(\Omega^{(0,\bullet)}(\rGr_{n,m}))) \simeq \Lambda_q(V) = \cT(V)/\langle S_q^2(V)\rangle.
  \]
  By definition, $\fB(V,-\sigma) = \cT(V) / \ker \mathfrak{S}^{-\sigma}$.
  By Corollary~\ref{cor:kerS2}, we have $\ker \mathfrak{S}^{-\sigma}_2 = S^2_q(V)$.
  Therefore
  $\langle S_q^2(V)\rangle$ is a sub-ideal in~$\ker \mathfrak{S}^{-\sigma}$,
  and hence there is an homogeneous ideal~$I$ in~$\Lambda_q(V)$ such that
  \begin{equation}\label{eq:BVLV}
    \fB(V,-\sigma) \simeq \Lambda_q (V) / I.
  \end{equation}
  In particular,  it follows that  $\fB(V,-\sigma)$ is finite-dimensional. 
  Let $p\colon \Lambda_q(V) \to \fB(V,-\sigma)$ be the canonical projection.
  Since both $\langle S_q^2(V)\rangle$ and $\ker \mathfrak{S}^{-\sigma}$ are
  $U_q(\fgl_m\oplus\fsl_{n-m})$-modules and graded ideals
  in~$\cT(V)$, it follows that $p$ is a graded $U_q(\fgl_m\oplus\fsl_{n-m})$-module algebra map.

  Let $d\in\Zee_{>0}$ be the largest integer such that $\fB_d(V,-\sigma) \neq 0$ and $\fB_{d+1}(V,-\sigma)=0$.
  It follows from~\eqref{eq:BVLV} that $d \leq N$.  For a~finite-dimensional Nichols algebra (and more generally, for
  a~finite-dimensional graded Hopf algebra in~$\rYD{H}$) Poincar\'e duality holds, which is to say,
  $\dim \fB_{d-k}(V,-\sigma) = \dim \fB_{k}(V,-\sigma)$,
  for   $k=0,\ldots,\lfloor d/2 \rfloor$, see~\cite[Proposition~3.2.2]{AndruskiewitschGrana1999}.
  In particular, $\dim \fB_0(V,-\sigma) = \dim \fB_d(V,-\sigma) = 1$. The quantum version of Howe
  duality~\eqref{eq:qHowe} gives us a decomposition of~$\Lambda_q(V)$ into a~sum of irreducible
  $U_q(\fgl_m\oplus\fsl_{n-m})$-modules. There are exactly two $1$-dimensional
  $U_q(\fgl_m\oplus\fsl_{n-m})$-modules in this decomposition, namely, $\Lambda_q^0(V)$ and
  $\Lambda_q^N(V)$. Therefore $p(\Lambda_q^0(V)) = \fB_0(V,-\sigma)$ and
  $p(\Lambda_q^N(V)) = \fB_d(V,-\sigma)$, implying that $d = N$.

  Denote by $\wedge$ the multiplication in~$\Lambda_q(V)$ and by $\wedge_\fB$ the multiplication in~$\fB(V,-\sigma)$.
  By~\cite[Proposition~4.11]{MTSUK}, we know that $\Lambda_q(V)$ is a Frobenius algebra. In particular  there exists a nondegenerate bilinear form
  \[
  (\cdot,\cdot)\colon \Lambda_q(V) \otimes \Lambda_q(V) \to \Cee,
  \]
  and an~element
  $\vol\in\Lambda_q^N(V)$ such that if $(v,v^c) = 1$ for $v,v^c\in\Lambda_q(V)$ then ${v\wedge v^c = \vol}$.

  Assume that $v\in I$ and $v\neq0$. Then 
  \[
    p(v\wedge v^c) = p(v) \wedge_\fB p(v^c) = 0.
  \]  
  On the other hand, since $\vol\in\Lambda_q^N(V)$,
  \[
    p(v\wedge v^c) = p(\vol) \neq 0,
  \]
  a contradiction. Hence there is no such~$v$, and $I = \langle 0 \rangle$ and
  $\Lambda_q(V) \simeq \fB(V,-\sigma)$.

  The proof that $\Phi(\Omega_q^{(\bullet,0)}(\rGr_{n,m}))$ is a Nichols algebra is analogous since $V^{(1,0)}$ and
  $V^{(0,1)}$ are dual $U_q(\fgl_m\oplus\fsl_{n-m})$-modules.
\end{proof}

\begin{cor}
  The calculi $\Omega_q^{(\bullet,0)}(\rGr_{n,m})$ and~$\Omega_q^{(0,\bullet)}(\rGr_{n,m})$ are Nichols algebras.
\end{cor}

\section{The general case of the irreducible quantum flag manifolds}\label{sec:other}

At present we do not have a quantum principal bundle description of the Hecken\-berger--Kolb calculi for irreducible quantum flag manifolds outside the $A$-series. However, we expect that such a description can be obtained in the same manner as for the quantum Grassmannians. Hence we
conjecture that the main results of the previous section hold for all irreducible quantum flag manifolds.

\begin{conj}[Analogue of Corollary~\ref{cor:VisYD}]\label{conj:YD}
  Let $\cO_q(G/L_S)$ be an irreducible quantum flag manifold and $\Omega^{(1,0)}_q(G/L_S)$,
  $\Omega^{(0,1)}_q(G/L_S)$ be the Heckenberger--Kolb first-order differential
  calculi. Then $\Phi(\Omega^{(1,0)}_q(G/L_S))$ and $\Phi(\Omega^{(0,1)}_q(G/L_S))$ 
  are Yetter--Drinfeld modules over~$\cO_q(L_S)$.
\end{conj}

\begin{conj}[Analogue of Theorem~\ref{th:mainNichols}]\label{conj:Nichols}
  For every irreducible quantum flag manifolds~$\cO_q(G/L_S)$ the maximal prolongations of the pair of Heckenberger--Kolb first-order differential calculi
  are Nichols algebras.
\end{conj}

In what follows we outline a strategy for proving Conjecture~\ref{conj:Nichols},  under the assumption that
Conjecture~\ref{conj:YD} is true. It is enough to calculate the eigenvalues of the braiding on the
irreducible components of the second tensor power of~$\Phi(\Omega_q^{(0,1)}(G/L_S))$, as discussed in  Lemma~\ref{lemma:YDbraidingOnHWV01}. This is
in contrast to the more involved proof of Theorem~\ref{th:mainNichols} for the quantum Grassmannians.

\subsection{Some general results for the irreducible quantum flag manifolds}
\subsubsection{The spectrum for $R$-matrices}
Let $V_\lambda$ be the irreducible type-1 representation of~$U_q(\fg)$ with the highest weight~$\lambda$.
Recall that by~\cite[\S8.4.3, Corollary~23]{KSLeabh}
the $R$-matrix braiding~\eqref{eq:RmatBraiding} acts on~$V_\mu$, an irreducible component with  highest
weight~$\mu$ in the decomposition of~$V_\lambda\otimes V_\lambda$, as multiplication by the scalar
\begin{equation}\label{eq:RmatEigen}
  \pm q^{-(2(\lambda,\lambda+2\rho)-(\mu,\mu+2\rho))/2},
\end{equation}
where $\rho$ is the half-sum of positive roots in~$\fg$.

\subsubsection{Quantum exterior algebras}
Let $\cO_q(G/L_S)$ be an irreducible quantum flag manifold and let $\Omega^{(0,1)}_q(G/L_S)$ and
$\Omega^{(1,0)}_q(G/L_S)$ be the corresponding pair of Hecken\-berger--Kolb first-order calculi. In what
follows we denote $V := \Phi(\Omega^{(0,1)}_q(G/L_S))$.

First recall that~$V$ is an irreducible $U_q(\fl_S)$-module with highest weight~$\lambda$, which
is the same as in the classical case (see the list in~\cite[\S6]{HK} and Table~3.2 on p.~27
in~\cite{BastonEastwood}).
Second, $V$ is a flat $U_q(\fl_S)$-module (for the complete list of flat modules see~\cite{Zwicknagl2009}) and
the quantum exterior algebra~$\Lambda_q(V)$ is isomorphic to~$\Phi(\Omega^{(\bullet,0)}_q(G/L_S))$. 
Third, all irreducible components in $V\otimes V$ have multiplicity 1, see~\cite[Table~5]{VinbergOnishchik}.
Hence together with the previous facts we can determine which sign occurs in~\eqref{eq:RmatEigen}.

We also present pictorial descriptions of the (quantum) irreducible flag manifolds. Note that we denote them by the same
symbols as in~\cite[Table~1]{Fano}, but the numbering of nodes in Dynkin diagrams follows~\cite[Table~1]{VinbergOnishchik}.
In what follows, $\fk_S$ denotes the maximal semisimple ideal of~$\fl_S$.

\subsubsection{Hecke-type braidings}

Let $V$ be a Yetter--Drinfeld module and $\sigma$ be the corresponding braiding.
The braiding $\sigma$ is of \emph{Hecke type} if
$(\sigma - \lambda)(\sigma +1) = 0$, for some non-zero scalar~$\lambda$. If $\lambda$ is not a root of unity
or if $\lambda = 1$, then the Nichols algebra is quadratic, see for example~\cite[Proposition~2.3]{Andruskiewitsch2004remarks}.
For example, the braiding~\eqref{eq:YDbraiding} for the quantum Grassmannians~$\cO_q(\mathrm{Gr}_{n,m})$ is of Hecke type
only for the special case of the quantum projective spaces (when $m=1$ or $m=n$).
In this case~$\Lambda_q^2(V)$ and $S^2_q(V)$ are irreducible $U_q(\fl_S)$-modules,
see~\S\ref{sec:CPn} for details.

\subsubsection{Quantum Grassmannians~$\cO_q(\mathrm{Gr}_{n,m})$} 
Consider the pictorial description of the (quantum) Levi subalgebra~$\fl_S$ corresponding to the crossed node
\[
\dynkin[%
  labels={1,,,,,,1},
  scale=2.5
]A{oo..oxo..oo}
\]
where in addition the numbered nodes determine the highest weight of the adjoint representation of~$\fg$.
In this case have that 
\[
  \fg=\fsl_{n},\quad \fl_S = \fgl_m\oplus \fsl_{n-m},\quad  \fk_S = \fsl_m\oplus\fsl_{m-n}, 
  \quad  V = V_{\varpi_1}\boxtimes V_{\varpi_1}^*.
\]
For $n>m>1$, the decomposition~\eqref{eq:qSLambda} of~$V\otimes V$ into quantum symmetric and antisymmetric parts with respect to diagonolised
$R$-matrix braiding is as follows
\[
  S_q^2(V) \simeq \left(V_{2\varpi_1}\boxtimes V_{\varpi_2}^*\right) \oplus \left(V_{\varpi_2}\boxtimes V_{2\varpi_1}^*\right),\quad
  \Lambda_q^2(V) \simeq  \left(V_{2\varpi_1}\boxtimes V_{2\varpi_1}^*\right) \oplus \left(V_{\varpi_2}\boxtimes V_{\varpi_2}^*\right).
\]
Conjecture~\ref{conj:YD} is proved in Corollary~\ref{cor:VisYD} and Conjecture~\ref{conj:Nichols} is proved in Theorem~\ref{th:mainNichols}.

\subsubsection{Quantum projective spaces~$\cO_q(\mathbb{CP}^n)$}\label{sec:CPn}
Consider the pictorial description of the (quantum) Levi subalgebra~$\fl_S$ corresponding to the crossed node
\[
  \dynkin[%
  labels={1,,,1},
  scale=2.5
]A{oo....ox}
\]
where in addition the numbered nodes determine the highest weight of the adjoint representation of~$\fg$.
In this case recall have 
\[
  \fg=\fsl_{n+1},\quad \fl_S = \fgl_n,\quad \fk_S = \fsl_n, \quad V = V_{\varpi_1}.
\]
For $n>1$, the decomposition~\eqref{eq:qSLambda} of~$V\otimes V$ into quantum symmetric and antisymmetric parts with respect to diagonolised
$R$-matrix braiding is as follows
\[
  S_q^2V \simeq V_{2\varpi_1},\qquad \Lambda_q^2V \simeq V_{\varpi_2}.
\]
Conjecture~\ref{conj:YD} is proved in Corollary~\ref{cor:VisYD} and Conjecture~\ref{conj:Nichols} is prooved
in Theorem~\ref{th:mainNichols}. Let us note that in this case the proof of Theorem~\ref{th:mainNichols} is
trivial since the Yetter--\/Drinfeld braiding is of Hecke type.
Moreover, the quasitriangular braiding gives the same Nichols algebra.

\subsubsection{Odd quantum quadrics~$\cO_q(\mathrm{Q}_{2n+1})$}
Consider the pictorial description of the (quantum) Levi subalgebra~$\fl_S$ corresponding to the crossed node
\[
\dynkin[%
  labels={,1,},
  scale=2.5
]B{xo....oo}
\]
where in addition the numbered node determines the highest weight of the adjoint representation of~$\fg$.
In this case we have that
\[
  \fg=\fo_{2n+1},\quad \fl_S = \fo_{2n-1}\oplus\fgl_1,\quad \fk_S = \fo_{2n-1},\quad V = V_{\varpi_1}.
  \]
For $n>2$, the decomposition~\eqref{eq:qSLambda} of~$V\otimes V$ into quantum symmetric and antisymmetric parts with respect to diagonalised
$R$-matrix braiding is as follows
\[
  S_q^2(V) \simeq  V_{2\varpi_1} \oplus V_0,\qquad  \Lambda_q^2(V) \simeq  V_{\varpi_2}.
\]

Assume that Conjecture~\ref{conj:YD} is true and let $\sigma$ be the corresponding Yetter--\/Drinfeld
braiding. If $\ker \mathfrak{S}^{-\sigma}_2 \simeq S_q^2V$, then, since $\Lambda_q^2(V)$ is irreducible, the
(rescaled) braiding~$-\sigma$ has two eigenvalues~$-1$ on~$S_q^2V$ and~$\lambda$ on~$\Lambda_q^2(V)$.
Therefore, the braiding~$-\sigma$ is of Hecke type and the
corresponding Nichols algebra is generated in degree~2 and isomorphic to~$\Lambda_qV$.
This implies Conjecture~\ref{conj:Nichols}.

Denote by $e_\mu$ the eigenvalue of the $R$-matrix braiding of~$U_q(\fk_S)$ on the irreducible
component with the highest weight~$\mu$ of~$V\otimes V$. We have
\[
  e_{2\varpi_1} = q,\qquad
  e_{0} = q^{-2(n-1)},\qquad
  e_{\varpi_2} = - q^{-1}.
\]

\subsubsection{Quantum Lagrangian Grassmannians~$\cO_q(\mathrm{L}_{n})$}
Consider the pictorial description of the (quantum) Levi subalgebra~$\fl_S$ corresponding to the crossed node
\[
\dynkin[%
  labels={2,},
  scale=2.5
]C{oo....ox}
\]
where in addition the numbered node determines the highest weight of the adjoint representation of~$\fg$.
In this case we have that
\[
  \fg=\fsp_{2n},\quad \fl_S = \fgl_{n},\quad \fk_S = \fsl_{n},\quad V = V_{2\varpi_1}.
\]
For $n>2$, the decomposition~\eqref{eq:qSLambda} of~$V\otimes V$ into quantum symmetric and antisymmetric parts with respect to diagonalised
$R$-matrix braiding is as follows
\[
  S_q^2V \simeq V_{4\varpi_1}\oplus V_{2\varpi_2}, \quad \Lambda_q^2V \simeq V_{2\varpi_1+\varpi_2}.
\]

Assume that Conjecture~\ref{conj:YD} is true and let $\sigma$ be the corresponding Yetter--\/Drinfeld
braiding. If $\ker \mathfrak{S}^{-\sigma}_2 \simeq S_q^2V$, then, since $\Lambda_q^2(V)$ is irreducible, the
(rescaled) braiding~$-\sigma$ has two eigenvalues~$-1$ on~$S_q^2V$ and~$\lambda$ on~$\Lambda_q^2(V)$.
Therefore, the braiding~$-\sigma$ is of Hecke type and the
corresponding Nichols algebra is generated in degree~2 and isomorphic to~$\Lambda_qV$.
This implies Conjecture~\ref{conj:Nichols}.

Denote by $e_\mu$ the eigenvalue of the $R$-matrix braiding of~$U_q(\fk_S)$ on the irreducible
component with the highest weight~$\mu$ of~$V\otimes V$. We have
\[
  e_{4\varpi_1} = q^{\tfrac{4}{n}(n-1)},\qquad
  e_{2\varpi_2} = q^{-\tfrac2n(n+2)},\quad
  e_{2\varpi_1+\varpi_2} = - q^{-\tfrac{4}{n}(n-1)}.
\]

\subsubsection{Even quantum quadrics~$\cO_q(\mathrm{Q}_{2n})$}
Consider the pictorial description of the (quantum) Levi subalgebra~$\fl_S$ corresponding to the crossed node
\[
\dynkin[%
  labels={,1,},
  scale=2.5
]D{xo....ooo}
\]
where in addition the numbered node determines the highest weight of the adjoint representation of~$\fg$.
In this case we have that
\[
  \fg=\fo_{2n},\quad \fl_S = \fo_{2(n-1)}\oplus\fgl_1,\quad  \fk_S = \fo_{2(n-1)},\quad  V = V_{\varpi_1}.
\]
For $n>2$, the decomposition~\eqref{eq:qSLambda} of~$V\otimes V$ into quantum symmetric and antisymmetric parts with respect to diagonalised
$R$-matrix braiding is as follows
\[
  S_q^2V \simeq V_{2\varpi_1} \oplus V_0,\qquad  \Lambda_q^2V \simeq V_{\varpi_2}.
\]
Assume that Conjecture~\ref{conj:YD} is true and let $\sigma$ be the corresponding Yetter--\/Drinfeld
braiding.  If $\ker \mathfrak{S}^{-\sigma}_2 \simeq S_q^2V$, then, since $\Lambda_q^2(V)$ is irreducible, the
(rescaled) braiding~$-\sigma$ has two eigenvalues~$-1$ on~$S_q^2V$ and~$\lambda$ on~$\Lambda_q^2(V)$.
Therefore, the braiding~$-\sigma$ is of Hecke type and the
corresponding Nichols algebra is generated in degree~2 and isomorphic to~$\Lambda_qV$.
This implies Conjecture~\ref{conj:Nichols}.

Denote by $e_\mu$ the eigenvalue of the $R$-matrix braiding of~$U_q(\fk_S)$ on the irreducible
component with the highest weight~$\mu$ of~$V\otimes V$. We have
\[
  e_{2\varpi_1} = q,\qquad
  e_{0} = q^{-2n+3},\qquad
  e_{\varpi_2} = - q^{-1}.
\]

\subsubsection{Quantum spinor varietes~$\cO_q(\mathrm{S}_{2n})$}
Consider the pictorial description of the (quantum) Levi subalgebra~$\fl_S$ corresponding to the crossed node
\[
\dynkin[%
  labels={,1,},
  scale=2.5
]D{oo....oxo}
\]
where in addition the numbered node determines the highest weight of the adjoint representation of~$\fg$.
In this case we have that
\[
  \fg=\fo_{2n},\quad \fl_S = \fgl_{n},\quad \fk_S = \fsl_{n},\quad V = V_{\varpi_2}.
\]
For $n>5$, the decomposition~\eqref{eq:qSLambda} of~$V\otimes V$ into quantum symmetric and antisymmetric parts with respect to diagonalised
$R$-matrix braiding is as follows
\[
  S_q^2V \simeq V_{2\varpi_2}\oplus V_{\varpi_4},\qquad  \Lambda_q^2V \simeq V_{\varpi_3+\varpi_1}.
\]
Assume that Conjecture~\ref{conj:YD} is true and let $\sigma$ be the corresponding Yetter--\/Drinfeld
braiding.  If $\ker \mathfrak{S}^{-\sigma}_2 \simeq S_q^2V$, then, since $\Lambda_q^2(V)$ is irreducible, the
(rescaled) braiding~$-\sigma$ has two eigenvalues~$-1$ on~$S_q^2V$ and~$\lambda$ on~$\Lambda_q^2(V)$.
Therefore, the braiding~$-\sigma$ is of Hecke type and the
corresponding Nichols algebra is generated in degree~2 and isomorphic to~$\Lambda_qV$.
This implies Conjecture~\ref{conj:Nichols}.

Denote by $e_\mu$ the eigenvalue of the $R$-matrix braiding of~$U_q(\fk_S)$ on the irreducible
component with the highest weight~$\mu$ of~$V\otimes V$. We have
\[
  e_{2\varpi_2} = q^{\tfrac{2}{n}(n-2)},\qquad
  e_{\varpi_4} = q^{-\tfrac{2}{n}(4n-5)},\qquad
  e_{\varpi_1+\varpi_3} = - q^{-\tfrac{2}{n}}.
\]

\subsubsection{Quantum Caley plane~$\cO_q(\mathbb{OP}^2)$}
Consider the pictorial description of the (quantum) Levi subalgebra~$\fl_S$ corresponding to the crossed node
\[
\dynkin[%
  labels={,1,},
  scale=2.5
]E{ooooox}
\]
where in addition the numbered node determines the highest weight of the adjoint representation of~$\fg$.
In this case we have that
\[
  \fg=\fe_{6},\quad \fl_S = \fo_{10}\oplus\fgl_1,\quad \fk_S = \fo_{10},\quad V = V_{\varpi_5}.
\]
The decomposition~\eqref{eq:qSLambda} of~$V\otimes V$ into quantum symmetric and antisymmetric parts with respect to diagonalised
$R$-matrix braiding is as follows
\[
  S_q^2V \simeq V_{2\varpi_5}\oplus V_{\varpi_1},\qquad  \Lambda_q^2V \simeq V_{\varpi_3}. 
\]
Assume that Conjecture~\ref{conj:YD} is true and let $\sigma$ be the corresponding Yetter--\/Drinfeld
braiding.  If $\ker \mathfrak{S}^{-\sigma}_2 \simeq S_q^2V$, then, since $\Lambda_q^2(V)$ is irreducible, the
(rescaled) braiding~$-\sigma$ has two eigenvalues~$-1$ on~$S_q^2V$ and~$\lambda$ on~$\Lambda_q^2(V)$.
Therefore, the braiding~$-\sigma$ is of Hecke type and the
corresponding Nichols algebra is generated in degree~2 and isomorphic to~$\Lambda_qV$.
This implies Conjecture~\ref{conj:Nichols}.

Denote by $e_\mu$ the eigenvalue of the $R$-matrix braiding of~$U_q(\fk_S)$ on the irreducible
component with the highest weight~$\mu$ of~$V\otimes V$. We have
\[
  e_{2\varpi_5} = q^{\tfrac{25}{5}},\qquad
  e_{\varpi_1} = q^{-\tfrac{27}{2}},\qquad
  e_{\varpi_3} = - q^{-\tfrac{3}{4}}.
\]

\subsubsection{Quantum Freudenthal variety~$\cO_q(\mathrm{F})$}
Consider the pictorial description of the (quantum) Levi subalgebra~$\fl_S$ corresponding to the crossed node
\[
\dynkin[%
  labels={1,,,,,,},
  scale=2.5
]E{oooooox}
\]
where in addition the numbered node determines the highest weight of the adjoint representation of~$\fg$.
In this case we have that
\[
  \fg=\fe_{7},\quad \fl_S = \fe_{6}\oplus\fgl_1,\quad \fk_S = \fe_{6}, \quad V = V_{\varpi_5}.
\]
The decomposition~\eqref{eq:qSLambda} of~$V\otimes V$ into quantum symmetric and antisymmetric parts with respect to diagonalised
$R$-matrix braiding is as follows
\[
  S_q^2V \simeq V_{2\varpi_5} \oplus V_{\varpi_1},\qquad \Lambda_q^2V \simeq  V_{\varpi_4}.
\]
Assume that Conjecture~\ref{conj:YD} is true and let $\sigma$ be the corresponding Yetter--\/Drinfeld
braiding.  If $\ker \mathfrak{S}^{-\sigma}_2 \simeq S_q^2V$, then, since $\Lambda_q^2(V)$ is irreducible, the
(rescaled) braiding~$-\sigma$ has two eigenvalues~$-1$ on~$S_q^2V$ and~$\lambda$ on~$\Lambda_q^2(V)$.
Therefore, the braiding~$-\sigma$ is of Hecke type and the
corresponding Nichols algebra is generated in degree~2 and isomorphic to~$\Lambda_qV$.
This implies Conjecture~\ref{conj:Nichols}.

Denote by $e_\mu$ the eigenvalue of the $R$-matrix braiding of~$U_q(\fk_S)$ on the irreducible
component with the highest weight~$\mu$ of~$V\otimes V$. We have
\[
  e_{2\varpi_5} = q^{\tfrac{2}{3}},\qquad
  e_{\varpi_1} = q^{-\tfrac{28}{3}},\qquad
  e_{\varpi_4} = - q^{-\tfrac{4}{3}}.
\]

\subsection{Nichols algebras and Weyl groupoids}

We are thankful to S.~Lentner who brought the following observation
to our attention. As shown in~\cite{AndruskiewitschHeckenbergerSchneider2010}, any Nichols algebra is controlled (but not
completely determined) by a Weyl groupoid.  Weyl groupoids are generalisations of Weyl groups motivated by
Serganova's work on generalised root systems of basic Lie superalgebras~\cite{Serganova1996} and examples
coming from from Nichols algebras, see~\cite{HeckenbergerYamane2008}, cf.~\cite{SergeevVeselov2011} and~\cite{GorelikHinichSerganova}.

A crystallographic arrangement~$\cA$ is a finite set of hyperplanes in $X:=\Ree^r$ which can be described
as kernels of a given set of root vectors~$\alpha_1,\ldots,\alpha_k\in V^\ast$ satisfying certain
properties (see~\cite{OrlikTerao1992,Cuntz2011} for an explicit definition).  Let $Y$ be a subspace of~$X$. Then
the restriction~$\cA^Y$ is defined to be the set of all hyperplanes of the form $Y\cap H$ for $H\in\cA$. In
general, $\cA^Y$ is not a crystallographic arrangement. When~$Y$ is an intersection of a subset of hyperplanes
of~$\cA$, $\cA^Y$ is called a parabolic restriction. In this case $\cA^Y$ is again a crystallographic arrangement. As was shown in~\cite{Cuntz2011}, we can associate a crystallographic arrangement to a given Weyl groupoid.

The set of simple roots~$I$ of the crystallographic arrangements of the Nichols
algebra~$\fB(V)$ enumerates the irreducible Yetter--Drinfeld submodules~$V_i$ of~$V=\bigoplus_{i\in I}V_i$.  Let
$J\subset I$ be a subset of simple roots for~$\fB(V)$ and denote by $V_J := \bigoplus_{i\in J} V_i$ the corresponding
Yetter--Drinfeld submodule of~$V$. Consider the associated sub-algebra of coinvariant elements 
$\fB(\tilde V) := \fB(V)^{\co(\fB(M_J))}$. The main result of~\cite{CuntzLentner2017} is that the root system
of~$\fB(\tilde V)$ is the parabolic restriction~$\cA^Y$, where $Y$ is the kernel of~$J$.
In particular, consider the Nichols algebra $\fB(V)$ associated
to the Borel part~$\mathfrak{u}_q(\fg)^+$ of Lusztig's small quantum group~$\mathfrak{u}_q(\fg)$ at a root of unity~\cite{Lusztig1990:L1,Lusztig1990:L2}. The restricted Nichols algebra~$\fB(\tilde V)$ belongs to the category of
Yetter--Drinfeld modules of~$u_q(\fg_J)$, where $\fg_J$ is the Lie subalgebra of~$\fg$ generated by the simple roots from~$J$.

The Nichols algebras corresponding to the Heckenberger--\/Kolb calculi for quantum Grassmannians have similar relations to these Nichols algebras. Thus we expect to see analogous behaviour for the Nichols algebras of the Heckenberger--Kolb calculi for all irreducible flag manifolds, assuming that Conjecture~\ref{conj:Nichols} is correct.

\providecommand{\bysame}{\leavevmode\hbox to3em{\hrulefill}\thinspace}
\providecommand{\MR}{\relax\ifhmode\unskip\space\fi MR }
\providecommand{\MRhref}[2]{%
  \href{http://www.ams.org/mathscinet-getitem?mr=#1}{#2}
}
\providecommand{\href}[2]{#2}

\end{document}